\documentclass[11pt,twoside]{amsart}
\usepackage[utf8]{inputenc}
\usepackage{amsmath}
\usepackage{amsfonts}
\usepackage{amssymb}
\usepackage{amsthm}
\usepackage[english]{babel}
\usepackage{hyperref}
\usepackage[alphabetic]{amsrefs}
\usepackage{geometry}
\usepackage{graphicx}
\usepackage{xcolor}
\usepackage{color}
\usepackage{mathrsfs}
\usepackage[indent]{parskip}
\usepackage{paralist}

\hypersetup{
	colorlinks   = true, 
	urlcolor     = blue, 
	linkcolor    = blue,
	citecolor    = blue 
}

\usepackage{mathtools}
\mathtoolsset{showonlyrefs}

\makeatletter
\def\thm@space@setup{%
	\thm@preskip=\parskip \thm@postskip=0pt
}
\makeatother
\usepackage{etoolbox}
\AtBeginEnvironment{proof}{\vspace*{-0.5\parskip}}

\usepackage{geometry}     % See geometry.pdf to learn the layout options. There are lots.
\usepackage{graphicx}
\usepackage{color}

\usepackage{amsthm}
\usepackage{amsmath}
\usepackage{amsfonts}
\usepackage{enumerate}
\usepackage{amssymb}
\usepackage{epstopdf}
\usepackage{hyperref}
\usepackage{comment}

\newtheorem{theorem}{Theorem}[section]
\newtheorem{proposition}[theorem]{Proposition}
\newtheorem{proposition*}{Proposition}
\newtheorem{lemma}[theorem]{Lemma}
\newtheorem{corollary}[theorem]{Corollary}

\newcommand{\eps}{\varepsilon}

\DeclareMathOperator{\spt}{spt}

\theoremstyle{definition}
\newtheorem{definition}[theorem]{Definition}

\newtheorem{remark}[theorem]{Remark}
\DeclareMathOperator*{\area}{area}
\newcommand{\dist}{\operatorname{dist}}
\newcommand{\II}{II}

\newcommand{\pr}{\partial}
\DeclareMathOperator*{\supp}{spt}

\theoremstyle{plain}
\newtheorem{thm}[theorem]{Theorem}
\newtheorem{lem}[theorem]{Lemma}
\newtheorem{pro}[theorem]{Proposition}

\newtheorem*{thm*}{Theorem}
\newtheorem*{lem*}{Lemma}
\newtheorem*{clm*}{Claim}

\theoremstyle{definition}
\newtheorem{defn}[theorem]{Definition}

\newtheorem{rmk}[theorem]{Remark}

\numberwithin{equation}{section}

\newcommand{\mf}{manifold}\newcommand{\mfs}{manifolds}

\newcommand{\Rm}{Riemannian}

\newcommand{\mt}{metric} 
\newcommand{\st}{such that}

\newcommand{\te}{there exist}\newcommand{\tes}{there exists} \newcommand{\Tes}{There exists}\newcommand{\fa}{for all}
\newcommand{\tf}{Therefore,} \newcommand{\hn}{Hence,}
\newcommand{\wrt}{with respect to}
\newcommand{\bbr}{\mathbb{R}}
\newcommand{\bbn}{\mathbb{N}}

\newcommand{\mres}{\scalebox{1.8}{$\llcorner$}}
\newcommand{\ra}{\rightarrow}
\newcommand{\fn}{function}

\newcommand{\sps}{Suppose}
\newcommand{\del}{\partial}
\newcommand{\seq}{sequence}\newcommand{\w}{with}
\newcommand{\cts}{continuous} 
\newcommand{\bF}{\mathbf{F}} 
\newcommand{\bM}{\mathbf{M}} 

\newcommand{\cm}{\mathcal{C}(M)}\newcommand{\ah}{\mathcal{A}^h}

\newcommand{\cA}{\mathcal{A}}\newcommand{\cC}{\mathcal{C}}\newcommand{\cF}{\mathcal{F}}\newcommand{\cH}{\mathcal{H}}\newcommand{\cS}{\mathcal{S}}\newcommand{\cV}{\mathcal{V}}

\newcommand{\sF}{\mathscr{F}}\newcommand{\sI}{\mathscr{I}}

\newcommand{\de}{\delta}\newcommand{\ve}{\varepsilon}\newcommand{\et}{\eta}\newcommand{\vp}{\varphi}\newcommand{\la}{\lambda}\newcommand{\rh}{\rho}\newcommand{\si}{\sigma}\newcommand{\tht}{\theta}\newcommand{\ta}{\tau}\newcommand{\ch}{\chi}\newcommand{\Ph}{\Phi}\newcommand{\Ps}{\Psi}\newcommand{\Om}{\Omega}\newcommand{\Si}{\Sigma}\newcommand{\Tht}{\Theta}\newcommand{\na}{\nabla}

\newcommand{\nm}[1]{\left\|#1\right\|}\newcommand{\md}[1]{\left|#1\right|}\newcommand{\db}[1]{[\![#1]\!]}
\newcommand{\vol}{\operatorname{Vol}}

\begin{document}
\title[On the existence of PMC hypersurfaces in non-compact manifolds]{On the existence of PMC hypersurfaces in complete manifolds of finite volume}

\author[G. Chambers]{Gregory R. Chambers}
\address{Department of Mathematics, Rice University, Houston, TX, USA}
\email{gchambers@rice.edu}

\author[A. Dey]{Akashdeep Dey}
\address{Chennai Mathematical Institute, Chennai, TN, India}
\email{akashdeep@cmi.ac.in}

\author[Y. Liokumovich]{Yevgeny Liokumovich}
\address{Department of Mathematics, University of Toronto, Toronto, ON, Canada}
\email{ylio@math.toronto.edu}

\author[X. Zhou]{Xin Zhou}
\address{Department of Mathematics, 531 Malott Hall, Cornell University, Ithaca, NY 14853, USA}
\email{xinzhou@cornell.edu}

\author[J. Zhu]{Jonathan J. Zhu}
\address{Department of Mathematics, University of Washington, Seattle, WA, USA}
\email{jonozhu@uw.edu}

\date{}

\maketitle

\begin{abstract}
Let $M^{n+1}$, $n+1\geq 3$, be a complete Riemannian manifold of finite volume, and suppose that $h\in C^{\infty}(M) \cap L^1(M)$ has the property that the width of $M$ \wrt\ the $\ah$ functional, $W^h(M)$, is positive.  Then there exists a hypersurface $\Si\subset M$ with prescribed mean curvature $h$, satisfying $0<\cH^n(\Si)<\infty$. Furthermore, $\Si$ is a closed subset of $M$ and it is smooth and almost embedded outside a singular set of Hausdorff dimension $\leq n-7$. The condition $W^h(M)>0$ is satisfied if $\|h\|_{L^1(M)}< W(M)$ (where $W(M)$ is the width of $M$ \wrt\ the area functional) or $\|h\|_{L^\infty(M)}< [\sI]_0$ (where $\sI$ is the isoperimetric profile of $M$ and $[\sI]_0:=\sup_{v>0} \frac{\sI(v)}{v}$).
\end{abstract}

\section{Introduction}

 The existence of minimal and prescribed mean curvature (PMC) hypersurfaces is a fundamental problem in geometric analysis, with deep connections to scalar curvature geometry, general relativity, and the calculus of variations. For closed Riemannian manifolds, the Almgren--Pitts min-max theory \cite{pit} provides a powerful framework for constructing minimal hypersurfaces.
Remarkable progress has been achieved for minimal hypersurfaces: 
in dimensions $3 \leq n+1 \leq 7$, 
every closed Riemannian manifold contains infinitely many closed embedded minimal hypersurfaces \cite{son}, and for 
a generic metric we have density and equidistribution results \cite{imn, mns}, as well as Morse index and multiplicity bounds \cite{mn-mult-one, zhou-mult}. The extension of min-max methods to hypersurfaces with prescribed mean curvature was achieved by Zhou--Zhu \cite{zz}, who proved that every closed manifold admits a smooth, closed, almost embedded hypersurface of any prescribed constant mean curvature.

The situation for non-compact manifolds is considerably more delicate. In \cite{cl}, Chambers and Liokumovich established that every complete non-compact manifold of finite volume contains a (possibly non-compact) minimal hypersurface of finite volume. Subsequently, Song \cite{son2} extended Yau's conjecture to complete manifolds that are ``thick at infinity,'' proving a dichotomy: either there are infinitely many saddle points of the $n$-volume functional, or there are none. Song's work also provided an alternative short proof of the existence of finite volume minimal hypersurfaces in finite volume manifolds.

For PMC hypersurfaces in non-compact settings, Mazurowski \cite{maz} developed techniques for applying one-parameter min-max theory when the prescription function $h$ is asymptotic to a positive constant near infinity, obtaining existence results in $\mathbb{R}^{n+1}$ and constant mean curvature surfaces in asymptotically flat $3$-manifolds. More recently, Stryker \cite{str} established a comprehensive min-max theory for PMC hypersurfaces in non-compact manifolds under the assumption that the prescription function $h$ does not change sign outside a compact set. 

Our main result also requires an assumption on the volume of the manifold; to this end we will say that $(M,g)$ has sublinear volume growth if
there exists $p \in M$
with $\lim\limits_{r \rightarrow \infty} \frac{\mathrm{Vol}(B_r(p))}{r} \allowbreak = 0$.
In particular, a manifold of finite volume has this property.
In this paper, we prove the existence of PMC hypersurfaces in complete Riemannian manifolds of sublinear volume growth under a set of hypotheses that do not require sign restrictions on $h$. Our main result is the following:

\begin{thm}\label{main theorem}
Let $M^{n+1}$, $n+1 \geq 3$, be a complete Riemannian manifold with sublinear volume growth, $h \in C^{\infty}(M) \cap L^1(M)$, and $W^h(M) > 0$. Then there exists a $h$-PMC hypersurface $\Sigma \subset M$ with optimal regularity (see Definition \ref{d.h.PMC.opt}), satisfying $0 < \mathcal{H}^n(\Sigma) < \infty$.
\end{thm}

Here $W^h(M)$ denotes the width of $M$ with respect to the $\mathcal{A}^h$ functional (see Definition~\ref{d.W^h(M)>0}). 

%In the case when $(M,g)$ has finite volume we obtain the following corollary.
\begin{corollary}
Let $(M^{n+1},g)$, $n+1 \geq 3$, be a complete Riemannian manifold with sublinear volume growth, $h \in C^{\infty}(M) \cap L^1(M)$. Suppose either of the following conditions hold:
\begin{enumerate}
    \item $\|h\|_{L^1(M)} < W(M)$, where $W(M)$ is the width with respect to the area functional (see Definition \ref{d.width}), or
    \item $(M,g)$ has finite volume, $h \in L^\infty(M)$, $\|h\|_{L^{\infty}(M)} < [\mathscr{I}]_0$, where $\mathscr{I}$ is the isoperimetric profile of $M$ and $[\mathscr{I}]_0 = \sup_{v > 0} \frac{\mathscr{I}(v)}{v}$.
\end{enumerate}
Then $W^h(M) > 0$ and \tes\ a $h$-PMC hypersurface $\Sigma \subset M$ with optimal regularity (see Definition \ref{d.h.PMC.opt}), satisfying $0 < \mathcal{H}^n(\Sigma) < \infty$.
\end{corollary}
The Corollary follows by Lemma~\ref{l.conditions.for.positive.W^h.1} and
Lemma~\ref{l.conditions.for.positive.W^h.2}. 

\subsection*{Statement on AI use}
AI has not been used during the preparation or writing of this article.
\subsection*{Acknowledgments}
A.D. was partially supported by an NSERC Discovery Grant during his postdoctoral fellowship at the University of Toronto, and is currently partially supported by an Infosys Foundation grant.
Y.L. acknowledges the support of the Simons Fellowship and NSERC discrovery grant.
X.Z. acknowledges the support by NSF grant DMS-2506717 and a grant from the Simons Foundation. 
JZ was supported in part by NSF grant DMS-2439945. 

\section{Min-max construction}

In this section, we revisit the min-max construction of PMC hypersurfaces in a \textit{closed} Riemannian manifold $M$. For an open set $U \subset M$, $\cC(U)$ denotes the space of Caccioppoli sets in $U$. For $E\in \cC(M)$, we have the following notations:
\begin{itemize}
    \item \[P(E,U)=\sup\left\{\int_{E\cap U}\operatorname{div}X\:d\cH^{n+1}:X\in\mathfrak{X}^1_c(U)\text{ and }\nm{X}_{C^0}\leq 1\right\},\]
    \item \[ V^h(E,U) = \int_{E\cap U}h \:d\cH^{n+1} \]
    where $h\in L^1(M)$,
    \item \[ V(E,U)=\cH^{n+1}(E\cap U), \]
    \item and
    \[\cA^h(E,U)=P(E,U)-V^h(E,U).\]
\end{itemize}
%$V^h(E,U) = \int_{E\cap U}h \:d\cH^{n+1}$ (where $h\in L^1(M)$), $V(E,U)=\cH^{n+1}(E\cap U)$, $\cA^h(E,U)=P(E,U)-V^h(E,U)$. 
Further, $P(E,M)$, $V^h(E,M)$, $V(E,M)$, $\ah(E,M)$ will be written simply as $P(E)$, $V^h(E)$, $V(E)$, $\ah(E)$. Without loss of generality, we can assume that in the definition of the $\ah$ functional, $V^h(M)\geq 0$. $\cF$ denotes the flat metric on $\cm$; the $\bF$ \mt\ on $\cm$\ is defined by
\[ \bF(E_1,E_2)=\cF(E_1,E_2)+\bF(|\del E_1|,|\del E_2|), \]
where $|S|$ denotes the varifold associated with $S$. Here we use $\cV_n(M)$ to denote the space of all $n$-dimensional varifolds in $M$. For the purpose of min-max construction of critical points of $\cA^h$, we also need the following notion: the $\mathcal{VC}$-space, denoted by $\mathcal{VC}(M)$, defined in \cite[Definition 1.3]{wang2024existenceminimalspheress3} (see also \cite[\S 2]{mazurowski2024infinitelyhalfvolumeconstantmean}), is the space of pairs $(V, E)$ consisting of a varifold $V\in \cV_n(M)$ and a Caccioppoli set $E\in \cC(M)$, such that there exists a sequence $\{E_k\}\subset \cC(M)$ satisfying $E_k\to E$ in $\cC(M)$ and $|\partial E_k|\to V$ in $\cV_n(M)$. The $\bF$ metric naturally extends to a metric $\mathscr{F}$ on $\mathcal{VC}(M)$ as defined in \cite[Definition 1.3]{wang2024existenceminimalspheress3}:
\[ \mathscr{F}((V_1, E_1), (V_2, E_2)) = \bF(V_1, V_2) + \cF(E_1, E_2).\]
The $\cA^h$-functional also naturally extends as
\[ \cA^h(V, E; U) = \|V\|(U) - V^h(E, U). \]
We refer to \cite[\S 1.1]{wang2024existenceminimalspheress3} and \cite[\S 2]{mazurowski2024infinitelyhalfvolumeconstantmean} for more details. 

\begin{definition}\label{d.h.PMC.opt}
    Let $M$ be an arbitrary Riemannian manifold, $h \in C^\infty(M)$. $\Si \subset M$ is called a \textit{$h$-PMC hypersurface with optimal regularity} if the following conditions are satisfied:
    \begin{itemize}
        \item[(i)] $\Si$ is a closed, $n$-rectifiable set;
        \item[(ii)] \tes\ a closed set $\mathcal{S} \subset M$ (which will be called the singular part of $\Si$) of Hausdorff dimension at most $n-7$, \st\ $\Si' := \Si \setminus \cS$ (which will be called the regular part of $\Si$) is a smooth, almost embedded hypersurface;
        \item[(iii)] the first variation of $|\Si|$ is in $L^\infty_{\text{loc}}(\|\Si\|);$
        \item[(iv)] the mean curvature of $\Si'$ at $x \in \{h=0\} \cap \Si'$ is zero; \tes\ an oriented smooth $n$-manifold $\tilde{\Si}$, a proper immersion $\iota : \tilde{\Si} \to \{h \neq 0\} \setminus \cS$ and a smooth unit normal $\nu$ on $\tilde{\Si}$ \st\ $\iota(\tilde{\Si}) = \{h \neq 0\} \cap \Si'$ and the mean curvature vector $H(x)$ (of the immersion $\iota$) at $x \in \tilde{\Si}$ is given by $H(x) = h(\iota(x)) \nu(x)$.
    \end{itemize}
\end{definition}

Let $\Om\subset M$ be an open set with smooth boundary, and $\sF(\Om)$ be the set of all $\cF$-continuous maps $\Ph:[0,1]\ra \cC(\Om)$ \st\ $\Ph(0)=\emptyset$ and $\Ph(1)=\Om$. The {\em $h$-PMC 1-width of $\Om$} is defined by
\[W^h(\Om)=\inf_{\Ph\in \sF(\Om)}\sup_{t\in [0,1]}\ah(\Ph(t),\Om). \]
Similarly, let $\widetilde{\sF}(\Om)$ (resp. $\widetilde{\sF}(\Om)^*$) be the set of all $\bF$-continuous (resp. $\cF$-continuous) maps $\Ph:[0,1]\ra \cC(M)$ \st\ $\Om \cap\Ph(0) =\emptyset$ and $\Om \setminus \Ph(1)=\emptyset$. The {\em reduced $h$-PMC 1-width of $\Om$} is defined by
\[\widetilde{W}^h(\Om)=\inf_{\Ph\in \widetilde{\sF}(\Om)}\sup_{t\in [0,1]}\ah(\Ph(t)).\]
By the interpolation results
\begin{equation}\label{e.tilde.W.*}
    \widetilde{W}^h(\Om)=\inf_{\Ph\in \widetilde{\sF}(\Om)^*}\sup_{t\in [0,1]}\ah(\Ph(t)).
\end{equation}
For $E\subset M$, let $E^c=M\setminus E$. For all $E\in \cm$,
\begin{align}
\ah(E,\Om)&=P(E,\Om)-V^h(E\cap\Om)\\
&\leq P(E)-V^h(E)+(V^h(E)-V^h(E\cap\Om))\\
&\leq \ah(E)+V^{|h|}(\Om^c).
\end{align}
Hence
\begin{equation}\label{e.W.tldW}
W^h(\Om)\leq \widetilde{W}^h(\Om)+V^{|h|}(\Om^c).
\end{equation}

A sequence $\{\Phi_i\}$ in $\widetilde{\sF}(\Om)$ is called a a {\em minimizing sequence} if $\sup_{t\in [0, 1]}\ah(\Phi_i(t)) \to \widetilde{W}^h(\Om)$. The associated {\em critical set} $\mathbf{C}(\{\Ph_i\})$ is:
\[\mathbf{C}(\{\Ph_i\}) = \left\{(V, E)\in \mathcal{VC}(M)\left|\,
\begin{aligned}
    & \exists \text{ a sequence $\{t_j\}\subset (0, 1)$ with $\ah(\Phi_{i_j}(t_j))\to \widetilde{W}^h(\Om)$,}\\
    & \text{such that } \mathscr{F}((|\partial \Phi_{i_j}(t_j)|, \Phi_{i_j}(t_j)), (V, E))\to 0.
\end{aligned}
\right\}\right..\]
In particular, we know that $\ah(V, E) = \widetilde{W}^h(\Om)$ if $(V, E) \in \mathbf{C}(\{\Ph_i\})$. A subsequence $\{\Phi_{i_j}(t_j)\}$ as above is usually called a {\em min-max sequence}.

\begin{pro}\label{p.pull.tight}
%There exists a critical sequence $S$ or minimizing sequence $\{\Ph_i\}$ in $\widetilde{\sF}(\Om)$ \st\ \fa\ $V\in \mathbf{C}(S) \text{ or } \mathbf{C}(\{\Ph_i\})$, either $V$ has $c$-bounded first variation (where $c=\sup |h|$) or $\nm{V}(\Om)=0$.
There exists a minimizing sequence $\{\Ph_i\}$ in $\widetilde{\sF}(\Om)$ such that for all $(V, E)\in \mathbf{C}(\{\Ph_i\})$, either $(V, E)$ is $\ah$-stationary (see \cite[Definition 1.1]{wang2024existenceminimalspheress3}), which implies that $V$ has $c$-bounded first variation (where $c=\sup |h|$), or $\nm{V}(\Om)=0$.
\end{pro}

\begin{proof}
The proof is almost identical to \cite[Proof of Proposition 5.1]{zz} or \cite[Theorem 2.7]{wang2024existenceminimalspheress3}. We modify the definition of $A^c_{\infty}$ in \cite[Section 5.1]{zz} and $A^L_0$ in \cite[(1.5)]{wang2024existenceminimalspheress3} as follows (see also \cite[Proof of Proposition 4.6]{mon}): First set $A^L = \{(V, E)\in \mathcal{VC}(M): \|V\|(M)\leq L\}$, and then let
\[
\begin{aligned}
    A^L_{\infty} & = \{(V, E) \in A^L: \de \ah_{V, E}(X) = \de V(X) - \int
    \md{\de V(X)}\leq c\int \md{X} d\mu_V,\, \forall X \in \mathfrak{X}(M)\} \\
    & \bigcup \{V\in A^L: \nm{V}(\Om)=0\}.
\end{aligned}
\]
%\[A^L_{\infty}=\{(V, E) \in A^L:\md{\de V(X)}\leq c\int \md{X} d\mu_V\text{, for any }X \in \mathfrak{X}(M)\}\bigcup \{V\in A^L: \nm{V}(\Om)=0\}.\]
Using this new definition, we can define the tightening map and obtain the pulled-tight minimizing sequence $\{\Ph_i\}$ exactly as in \cite[Proof of Proposition 5.1]{zz} or \cite[Theorem 2.7]{wang2024existenceminimalspheress3}. More precisely, we choose an arbitrary minimizing sequence $\{\Ph_i^*\}\subset \widetilde{\sF}(\Om)$ and define $\Ph_i(x)=\Ps_{\md{\del \Ph_i^*(x)}}(1)(\Ph_i^*(x))$, where $\Ps$ is the tightening map.
\end{proof}

\begin{defn}\label{d.prop.P.eta}
Let $\Om\subset M$ be an open set with smooth boundary. $\Om$ is said to satisfy property $(P_{\et})$ (for some $\et >0$) \wrt\ the $\cA^h$ functional if the following holds: 
\begin{itemize}
    \item[($P_{\et}$)] there exists $\Tht \Subset \Om$ \st\ for all $\Ph\in \widetilde{\sF}(\Om)$, there exists $t_0\in [0,1]$ so that $\ah(\Ph(t_0))\geq \widetilde{W}^h(\Om)$ and $\|\del \Ph(t_0)\|(\overline{\Tht})\geq \et$.
\end{itemize} 
\end{defn}

\begin{defn}\label{d.mdf.crit.set}
Suppose $\Om \subset M$ satisfies property $(P_{\et})$ \wrt\ the $\cA^h$ functional and $\Tht$ is as in Definition \ref{d.prop.P.eta}. The {\em modified critical set} $\mathbf{C}^*(\{\Ph_i\})$ of a minimizing sequence $\{\Ph_i\}\subset \widetilde{\sF}(\Om)$ is the set of all $(V, E)\in \mathcal{VC}(M)$ satisfying the following condition: there exists a min-max subsequence $\{\Ph_{i_j}(t_j)\}$, such that
\[
\mathscr{F}\big((|\partial \Phi_{i_j}(t_j)|, \Phi_{i_j}(t_j)), (V, E)\big) \to 0,\; \lim_{j\to \infty} \ah(\Ph_{i_j}(t_j))= \widetilde{W}^h(\Om),\; \|\del \Ph_{i_j}(t_j)\|(\overline{\Tht})\geq \et, \;\forall j.
\]
Note that Property $(P_{\et})$ ensures that $\mathbf{C}^*(\{\Ph_i\})$ is always nonempty.
%\textcolor{red}{(Modify.)} of all $V\in \cV_n(M)$ such that the following holds: there exists a subsequence $\{\Ph_{i_j}\}\subset \{\Ph_i\}$ \st\ 
%$$V= \lim_{j\to \infty}\md{\del \Ph_{i_j}(x_j)},\;   \lim_{j\to \infty} \cA^h(\Ph_{i_j}(x_j))= \widetilde{W}^h(\Om),\; \|\del \Ph_{i_j}(x_j)\|(\overline{\Tht})\geq \et \;\forall j.$$
\end{defn}

We now recall the notion of almost minimizing pairs; see \cite[Definition 3.2]{wang2024existenceminimalspheress3}, \cite[Definition 4.5]{mazurowski2024infinitelyhalfvolumeconstantmean}. Let $\nu$ be the $\mathcal F$ or $\mathbf F$ norms on $\cC(M)$. Fix an open subset $U\subset M$. Given $\eps, \de>0$, an $E\in \cC(M)$ is called {\em $(\ah, \eps, \de, \nu)$-almost minimizing} in $U$, if $E=E_0, E_1, \cdots, E_m$ is a sequence in $\cC(M)$, satisfying:
\begin{itemize}
    \item $\spt(E_i-E) \subset U$,
    \item $\nu(E_i, E_{i+1})<\de$,
    \item $\ah(E_i) \leq \ah(E)+\de$, 
\end{itemize}
then $\ah(E_m) \geq \ah(E) - \eps$. 

\begin{defn}
    A pair $(V, E)\in \mathcal{VC}(M)$ is called {\em $\ah$-almost minimizing with respect to a sequence $\{E_j\}\subset \cC(M)$ in $U$}, if there exist $\eps_j \to 0$ and $\de_j\to 0$, such that
\begin{itemize}
    \item $\mathscr{F}\big(E_j, (V, E)\big) \to 0$ as $j\to \infty$, and
    \item $E_j$ is $(\ah, \eps_j, \de_j, \mathcal{F})$-almost minimizing in $U$ for each $j$. 
\end{itemize}
\end{defn}

Recall that given $L\in \mathbb{N}$ and $p\in M$, a collection of annuli centered at $p$:
\[ \mathscr{C} = \{ An(p, s_1, r_1), \cdots, An(p, s_L, r_L)\}\]
is said to be {\em $L$-admissible} if $2r_{j+1} < s_j$ for all $1\leq j\leq L-1$; see for instance \cite[Definition 3.9]{wang2024existenceminimalspheress3}.

\begin{defn}
    Let $\mathscr{C}$ be a collection of concentric annuli. A pair $(V, E)\in \mathcal{VC}(M)$ is called {\em $\ah$-almost minimizing in $\mathscr{C}$ with respect to a sequence $\{E_j\}\subset \cC(M)$}, if there exist $\eps_j \to 0$ and $\de_j\to 0$, such that
\begin{itemize}
    \item $\mathscr{F}\big(E_j, (V, E)\big) \to 0$ as $j\to \infty$, and
    \item for each $j$, $E_j$ is $(\ah, \eps_j, \de_j, \mathcal{F})$-almost minimizing in at least one annulus in $\mathscr{C}$. 
\end{itemize}
\end{defn}

The following is a standard annulus-picking argument (see for instance \cite[Lemma A.3]{Colding-Gabai_Ketover18}). We note that \cite{wang2024existenceminimalspheress3} proved a significant generalisation involving countably many collections $\mathcal{P}_i$ (see \cite[Proposition D.1]{wang2024existenceminimalspheress3}).
\begin{lemma}
    \label{lem:an-pick-easy}
    Let $\mathcal{P}$ be a collection of annuli and suppose that for any $L$-admissible collection $\mathscr{C}$ of annuli centred at $p$, at least one annulus in $\mathscr{C}$ belongs to $\mathcal{P}$. Then there exists $r_p>0$ such that for all $s<r<r_p$, the annulus $An(p,s,r)$ is in $\mathcal{P}$.
\end{lemma}

%Given a pair $(V, E)\in \mathcal{VC}(M)$, and a sequence $\{E_j\}\subset \cC(M)$

By adapting the combinatorial arguments as in \cite{pit} and \cite{zz.cmc}, we have the following.

\begin{pro}\label{p.a.m.}
There exists an integer $L$ (depending only on the dimension $n$ and the number $p$ of min-max parameters) with the following significance: 

Suppose $\Om \subset M$ satisfies property $(P_{\et})$ \wrt\ the $\cA^h$ functional and $\{\Ph_i\}\subset \widetilde{\sF}(\Om)$ is a minimizing sequence. 
Then there is a min-max sequence $\{E_j = \Phi_{i_j}(t_j): t_j\in [0, 1]\}$ converging to some $(V, E)\in \mathbf{C}^*(\{\Ph_i\})$, such that
\begin{itemize}
    \item[(a1)] $\|V\|(\overline{\Tht}) \geq \et$, and
    \item[(a2)] $(V, E)$ is $\ah$-almost minimizing in any $L$-admissible family of annuli with respect to $\{E_j\}$.
\end{itemize}
%Consequently, there exists $(V, E)\in \mathbf{C}^*(\{\Ph_i\})$ such that
In particular, $(V,E)$ satisfies: %\textcolor{red}{(JZ: I don't think we make a new choice of $(V,E)$? feel free to switch back)}
\begin{itemize}
    \item[(b1)] $\|V\|(\overline{\Tht})\geq \et$, and
    \item[(b2)] $(V, E)$ is $\ah$-almost minimizing in small annuli, that is: for any $p\in M$, there exists $r_{am}(p)>0$, such that for any annulus $An(p, s, r)$ with $s<r<r_{am}(p)$, $(V, E)$ is $\ah$-almost minimizing in $An(p, s, r)$ with respect to some sequence $\{\Phi_{i_j}(t_j): t_j\in [0, 1]\}$.
\end{itemize}

\end{pro}
\begin{proof}
    We first prove (a1) and (a2). Suppose that these conclusions do not hold, then the arguments in \cite[Proof of Lemma 3.1]{wang2024existenceminimalspheress3} may be adapted as follows: We can find some $\epsilon_0>0$, such that for any $\delta>0$ and $i>1/\delta$, and any $t\in (0, 1)$ satisfying 
    \[ \|\partial\Phi_i(t)\|(\overline{\Tht}) \geq \eta,\quad \text{and} \quad \ah(\Phi_i(t))\geq  \widetilde{W}^h(\Om) - \epsilon_0, \]
    there exist an $L$-admissible collection $\mathscr{C}_{i, t}$ of annuli, such that $\Phi_i(t)$ is not $(\ah, \epsilon, \delta)$-almost minimizing in any annulus in $\mathscr{C}_{i, t}$. We can then use the combinatorial arguments as in \cite[Proof of Theorem 4.10]{pit}, \cite[Proof of Theorem 5.6]{zz.cmc} and \cite[Proof of Theorem 1.7]{zhou-mult} to obtain $\tilde{\Ph}\in \widetilde{\sF}(\Om)$ \st\ 
    \begin{equation}\label{e.pf.p.am}
    \sup\left\{\ah(\tilde{\Ph}(t)):t\in [0,1]\text{ such that } \|\del \tilde{\Ph}(t)\|(\overline{\Tht})\geq \et \right\} < \widetilde{W}^h(\Om).
    \end{equation}
    However, \eqref{e.pf.p.am} contradicts the assumption that $\Om $ satisfies property $(P_{\et})$ \wrt\ the $\cA^h$ functional.

    Finally, note that item (b1) trivially follows from item (a1). The implication from (a2) to (b2) follows from a standard argument as in \cite[Proof of Theorem 3.8]{wang2024existenceminimalspheress3} using \cite[Proposition D.1]{wang2024existenceminimalspheress3}. 

\end{proof}

%\begin{proof}
%We assume by contradiction that \tes\ a minimizing \seq\ $\{\Ph_i\}$ \st\ \fa\ $V\in \mathbf{C}^*(\{\Ph_i\})$ \tes\ $p\in M$ \st\ $V$ is not $h$-almost minimizing in small annuli centered at $p$. We also note that 	$\mathbf{C}^*(\{\Ph_i\})$ is a compact subset of $\cV_n(M)$. We can argue exactly as in \cite[Proof of Theorem 4.10]{pit} and \cite[Proof of Theorem 5.6]{zz.cmc} and obtain $\tilde{\Ph}\in \widetilde{\sF}(\Om)$ \st\ 
%\begin{equation}\label{e.pf.p.am}
%\sup\left\{\ah(\tilde{\Ph}(t)):t\in [0,1]\text{ such that } \|\del \tilde{\Ph}(t)\|(\overline{\Tht})\geq \et \right\} < \widetilde{W}^h(\Om).
%\end{equation} 
%(We need to replace the critical set by the modified critical set in that argument; $\tilde{\Ph}$ is obtained by deforming $\Ph_i$ for sufficiently large $i$.) However, \eqref{e.pf.p.am} contradicts the assumption that $\Om $ satisfies property $(P_{\et})$ \wrt\ the $\cA^h$ functional.
%\end{proof}

We summarize the results of Propositions \ref{p.pull.tight} and \ref{p.a.m.} by the following theorem:

\begin{thm}\label{t.a.m.}
Suppose $\Om \subset M$ satisfies property $(P_{\et})$ \wrt\ the $\cA^h$ functional. Then \tes\ $(V, E)\in \mathcal{VC}(M)$ \st
\begin{enumerate}[(i)]
	\item $(V, E)\in \mathbf{C}(\{\Ph_i\})$ for some minimizing sequence $\{\Ph_i\}$ in $\widetilde{\sF}(\Om)$, and hence $\ah(V, E) = \widetilde{W}^h(\Om)$;
    %$V\in \mathbf{C}(S) \text{ or } \mathbf{C}(\{\Ph_i\})$ for some critical sequence $S$ or minimizing sequence $\{\Ph_i\}$ in $\widetilde{\sF}(\Om)$;
	\item $(V, E)$ is $\ah$-stationary, and hence $V$ has $c$-bounded first variation (where $c=\sup |h|$);
    %$V$ has $c$-bounded first variation (where $c=\sup |h|$);

	\item $\nm{V}(\overline{\Tht})\geq \et $;
	\item $(V, E)$ is $\ah$-almost minimizing in small annuli. 

\end{enumerate}
\end{thm}

We now use approximation methods and compactness results for stable PMC hypersurfaces to establish the regularity of our min-max PMC hypersurfaces. For a discussion of stability and consequences for PMC hypersurfaces (particularly the subtly different notions of stability), the reader may consult Section \ref{sec:stability}. 

\begin{thm}[Regularity Theorem]\label{t.regularity}
    Let $h$ be a smooth function on a closed $(n+1)$-dimensional Riemannian manifold $M$. Suppose that $n\geq 2$, and let $(V, E)$ be as in Theorem \ref{t.a.m.}. Then $V \mres \{h \neq 0\} = |\del E|\mres \{h \neq 0\}$ and the support of $\|V\|$ is a $h$-PMC hypersurface with optimal regularity.
\end{thm}
\begin{proof}
    {\bf Part 1:} We first explain how to generalize the regularity result in \cite{zz} for generic $h\in C^\infty(M)$ to dimensions $n+1\geq 8$. Let $h$ be a generic smooth function given by \cite[Proposition 0.2]{zz}. We now outline the key modifications when $n+1\geq 8$.

    The first major difference in higher dimensions $n+1\geq 8$ is the regularity of the constrained minimizer \cite[Lemma 6.7]{zz}. Recall that we use constrained minimizers to construct {\em $h$-replacements} for $(V, E)$ in annuli $An(p, s, r)$ in which $(V,E)$ is $\ah$-almost-minimizing.  
    
    In fact, the constrained minimizer $E_j^*\in \mathcal{C}(M)$ associated with $E_j$ (here $\{E_j\}$ is the min-max sequence given in (a2) of Proposition \ref{p.a.m.}) is locally $\ah$-minimizing in any annuli $An(p, s, r)$ where (a2) holds true. By the classical regularity theory for PMC hypersurfaces (see for instance \cite{Mo03} or \cite[Section 2.2]{zz}), $\partial E_j^*$ is a properly embedded $h$-PMC hypersurface with a singular set of dimension at most $n-7$ in $An(p, s, r)$. Moreover, by the construction $\pr E^*_j$ will be $\ah$-stable in $An(p, s, r)$ (with respect to ambient deformations, hence as an immersion since $\pr E^*_j$ is properly embedded, away from its singular set). 
    
    Next, letting $j\to \infty$, after taking a subsequence, $E_j^* \to (V^*, E^*)$ under the $\mathscr{F}$-metric. Moreover, applying the stable compactness Theorem \ref{thm:pmc-compactness} in $An(p, s, r)$, we get that $\partial E_j^*$ converges locally smoothly, away from a singular set of dimension at most $n-7$, to an almost embedded $h$-PMC hypersurface $\Sigma^*$. Moreover, $\Sigma^*$ is $\ah$-stable (relative to ambient functions) in $An(p,s,r)$ (away from its singular set). In particular, we know that $V^*\llcorner An(p, s, r)$ is induced by $\partial E^*\llcorner An(p, s, r)$, and their supports coincide with $\Sigma^*$. The pair $(V^*, E^*)$ is the {\em $h$-replacement} for $(V, E)$ in $An(p, s, r)$. 

    Finally, using the existence and stability of the $h$-replacements, we can extend the procedure in \cite[Section 7]{zz} in the same way that Schoen-Simon \cite[Section 7]{schoen1981regularity} extends the minimal hypersurface min-max regularity to higher dimensions, which yields the desired regularity for $(V, E)$.

    In preparation for dealing with arbitrary functions $h$ in the following part, we highlight some stability properties of the min-max pair $(V,E)$. By an observation in \cite[\S 6]{Zhou15}, we know that the $h$-replacement $(V^*, E^*)$ coincide with $(V, E)$ inside $An(p, s, r)$ when $r$ is less than a uniform constant\footnote{For instance, $r_0$ can be chosen to be less than the injectivity radius of $M$, and so that $\partial B(p, r_0)$ has mean curvature greater than $c=\sup|h|$.} $r_0$ depending only on $(M,g)$ and $c=\sup|h|$. Hence we know that $(V, E)$ satisfies:
    \begin{itemize}
        \item[(S)] if $(V, E)$ is $\ah$-almost minimizing in any annulus $An(p, s, r)$ for $0<s<r<r_0$, then $(V, E)$ satisfies the desired regularity, and the corresponding $h$-PMC hypersurface $\Sigma$ is $\ah$-stable (relative to ambient functions) in $An(p, s, r)$ (away from its singular set).  
    \end{itemize}

    \vspace{0.5em}
    {\bf Part 2:} We now explain how to obtain the regularity result for a general $h\in C^\infty(M)$ by approximations. Pick a sequence of generic $h_k\in C^\infty(M)$ which converge smoothly to $h$, and consider the sequence of min-max problems associated with $\mathcal{A}^{h_k}$ as in Part 1. The min-max values $\widetilde{W}^{h_k}(\Om)$ will converge to $\widetilde{W}^h(\Om)$. Let $(V_k, E_k)$ be the associated min-max pair for each $\widetilde{W}^{h_k}(\Om)$. We may assume that $(V_k, E_k)$ converges to a pair $(V, E)$ under the metric $\mathscr{F}$, and hence $\ah(V, E) = \lim_{k\to\infty} \mathcal{A}^{h_k}(V_k, E_k) = \widetilde{W}^h(\Om)$.
    
    Consider $p\in M$. By Proposition \ref{p.a.m.}, we know that for any $L$-admissible family $\mathscr{C}$ of annuli centred at $p$, each $(V_k, E_k)$ is $\mathcal{A}^{h_k}$-almost minimizing in $\mathscr{C}$. The key is that $L$ and the constant $r_0$ in property (S) do not depend on $k$. 
    
     % Method 1: (\textcolor{magenta}{this uses only the simplified Lemma \ref{lem:an-pick-easy} to pick annuli after the limit}) 
    
    Take an $L$-admissible family $\mathscr{C}$ of annuli centred at $p$. Then up to choosing a subsequence, there is an annulus $An(p,s,r)$ in $\mathscr{C}$ so that each $(V_k, E_k)$ is $\mathcal{A}^{h_k}$-almost minimising in $An(p,s,r)$ (away from its singular set). By property (S), each $(V_k, E_k)$ satisfies the desired regularity, and the corresponding $h_k$-PMC hypersurface $\Sigma_k$ is $\mathcal{A}^{h_k}$-stable (relative to ambient functions) in $An(p, s, r)$. By the stable compactness Theorem \ref{thm:pmc-compactness}, a subsequence of the $\Sigma_k$ converges locally smoothly in $A(p,s,r)$ to a properly embedded $h$-PMC hypersurface $\Sigma$ away from a singular set of dimension at most $n-7$. We conclude that for any $L$-admissible family $\mathscr{C}$ of annuli centred at $p$, the limit $(V,E)$ has the desired regularity in some annulus in $\mathscr{C}$. 

    But now by Lemma \ref{lem:an-pick-easy}, there exists $r(p)>0$ such that $(V,E)$ has the desired regularity in every annulus $An(p, s, r)$ with $0<s<r<r(p)$. (One could also have argued using the full power of \cite[Proposition D.1]{wang2024existenceminimalspheress3} to pick annuli before taking a limit.)
    
    % Method 2: (\textcolor{magenta}{this uses the full power of \cite[Proposition D.1]{wang2024existenceminimalspheress3} to pick annuli before taking a limit. did not clean fully}) By a standard argument as in \cite[Proposition 5.1]{wang2024existenceminimalspheress3} using \cite[Proposition D.1]{wang2024existenceminimalspheress3}, we know that, up to a subsequence, the following holds true:
    % \begin{itemize}
    %     \item[(AM)] for each $p\in M$, there exists $r(p)>0$, such that for each annulus $An(p, s, r)$ with $0<s<r<r(p)$, $(V_k, E_k)$ is $\mathcal{A}^{h_k}$-almost minimizing in $An(p, s, r)$ for sufficiently large $k$.
    % \end{itemize}
    % By the above Property (S), we know that $(V_k, E_k)$ satisfy the desired regularity, and the corresponding $h_k$-PMC hypersurface $\Sigma_k$ is $\mathcal{A}^{h_k}$-stable in $An(p, s, r)$ as an immersion for $k$ sufficiently large. 

    % Finally, in each $An(p, s, r)$, letting $k\to\infty$, we know by the stable compactness Theorem \ref{thm:pmc-compactness} that $\Sigma_k$, up to a subsequence, converges locally smoothly to an \textcolor{blue}{almost} embedded $h$-PMC hypersurface $\Sigma$ outside a singular set of dimensional at most $n-7$. This establishes the regularity of $(V, E)$ in $An(p, s, r)$. 
    
    By taking a finite cover, this establishes the desired regularity, possibly except at finitely many discrete points $p_1,\cdots, p_m$. For $n\geq 7$, we are done as the latter set $\{p_i\}$ has dimension $0 \leq n-7$. For $2\leq n \leq 6$, the Removable Singularity Theorem \cite[Theorem B.1]{zhou-mult} completes the regularity of $(V,E)$ near the points $p_i$.
\end{proof}

From the proof of the above Theorem \ref{t.regularity}, we can observe: 

\begin{corollary}
    \label{cor:regularity-stable}
    Let $(M,g)$ be a closed Riemannian manifold, $h\in C^\infty(M)$, and let $\Si$ be as in Theorem \ref{t.regularity}. Define $\rh: M \to (0,1]$ by $\rh(p) = \min \{1, \text{inj}_p(M)\}$. Then, it follows from the above proof that for all $p \in M$, \tes\ $0< \si_p \le \rh(p)$, depending only on the restrictions of $g$ and $h$ to $B_{\rh(p)}(p)$, such that:
    \begin{enumerate}[(a)]
        \item $\Si$ is stable relative to ambient functions (away from its singular set) in $B_{\si_p}(p) \setminus \{p\}$;
        \item when $2 \le n \le 6$, $\Si$ is stable under ambient deformations in $B_{\si_p}(p) \setminus \{p\}$.
    \end{enumerate}
\end{corollary}

This stability will be later be used to take a further limit; indeed, we record:

\begin{remark}\label{r.lim.min-max.hyp}
    
By the compactness theorems discussed in the next section and the Removable Singularity Theorem in \cite[Theorem B.1]{zhou-mult}, it follows that: If  $\{V_i\}$ is a sequence of integral $n$-varifolds in $M$ such that $\|V_i\|(M)$ is uniformly bounded, $\Si_i = \text{spt} \|V_i\|$ is a $h$-PMC hypersurface satisfying the conclusions of Corollary \ref{cor:regularity-stable}, and $V_i \to V_\infty$ in the sense of varifolds, then $\text{spt}\|V_\infty\|$ is a $h$-PMC hypersurface with optimal regularity.
\end{remark}

\subsection{Preliminaries on stability}
\label{sec:stability}

Recall that the index form of a 2-sided $h$-PMC hypersurface $\Sigma^n$ in $M^{n+1}\hookrightarrow \mathbb{R}^N$ is 
\begin{equation}
\label{eq:2var}
\II_\Sigma(\varphi,\varphi) =\int_{\Sigma}\left( |\nabla\varphi|^2-\left(Ric^M(\nu, \nu)+|A^\Sigma|^2 + \partial_\nu h\right)\varphi^2\right)d\mu_{\Sigma}. 
\end{equation}

Consider an open set $U\subset M$. Stability of $\Sigma$ (in $U$) may often be defined by nonnegativity of the index form, but we take care to emphasise the admissible test functions $\varphi$. Denote the immersion $\iota:\Sigma\looparrowright M$. 
\begin{definition}
We say that $\Sigma$ is stable under ambient deformations (in $U$) if $\II_\Sigma(\varphi,\varphi)\geq 0$ for all test functions of the form $\varphi=\langle X,\nu\rangle$, where $X$ is a $C^1$ ambient vector field tangent to $M$, with compact support in $U$. 

We say that $\Sigma$ is stable as an immersion (in $U$) if $\II_\Sigma(\varphi,\varphi)\geq 0$ for all $C^1$ test functions $\varphi$ with compact support in $\iota^{-1}(U)$. 

Finally, $\Sigma$ is stable relative to ambient functions (in $U$) if $\II_\Sigma(\varphi,\varphi)\geq 0$ for all test functions $\varphi = \phi\circ\iota$, where $\phi$ is an ambient $C^1$ function with compact support in $U$.
\end{definition}

\begin{remark}
\label{rmk:2var-cts}
The expression (\ref{eq:2var}) is clearly continuous in smooth variations of $M, h, \varphi$ and $\Sigma$ (as an immersion). By standard approximation arguments, one may equivalently require Lipschitz regularity for the test functions above, instead of $C^1$. We also note that stability as an immersion certainly implies stability relative to ambient functions.
\end{remark}

The main compactness tool we use is the following consequence of the regularity theory of Bellettini-Wickramasekera \cite{bw}:
\begin{theorem}[\cite{bw}]
\label{thm:pmc-compactness}
    Let $n\geq2$ and consider an orientable Riemannian manifold $(U^{n+1},g)$ and function $h_k,h \in C^\infty(U)$, with $h_k\to h$ in $C^\infty_{loc}$. 

    Suppose that $V_k$ is a sequence of integral $n$-varifolds in $U$ and $\mathcal{S}_k$ is a sequence of closed sets of Hausdorff dimension at most $n-7$ such that $\Sigma_k = \spt V_k \setminus \mathcal{S}_k$ is a two-sided almost embedded $h_k$-PMC hypersurface, which is stable relative to ambient functions in $U\setminus \mathcal{S}_k$. Also assume that $\limsup_{k\to\infty} \|V_k\|(K) <\infty$ for each $K\subset\subset U$. 

    Then, possibly after passing to a subsequence, $V_k$ converges to an integral $n$-varifold $V$ in $U$ such that $\Sigma= \spt V \setminus \mathcal{S}$ is an almost embedded $h$-PMC hypersurface, where $\mathcal{S}$ is a closed set of Hausdorff dimension at most $n-7$. Moreover, $\Sigma$ is stable relative to ambient functions in $U\setminus\mathcal{S}$. 
    
    %\textcolor{magenta}{say something about the convergence?}
\end{theorem}

In the remainder of this subsection, we describe how one may obtain compactness results in low dimensions $2\leq n \leq 6$ by more elementary methods. We begin with the following sheeted convergence lemma under bounded curvature assumptions.

\begin{lemma}
\label{lem:bdd-cpctness}
Consider a sequence of Riemannian manifolds $(U, g_k)$ converging locally smoothly to a Riemannian manifold $(U,g)$, and $h_k\in C^\infty(U)$ converging locally smoothly to some $h\in C^\infty(U)$. 

 Suppose $\Sigma_k$ is a sequence of smooth, almost embedded, two-sided, $h_k$-PMC hypersurfaces in $(U,g_k)$ with $\pr \Sigma_k\cap U=\emptyset$. If the area of $\Sigma_k$ and the curvature $|A^{\Sigma_k}|$ are uniformly bounded on compact subsets, then up to a subsequence, $\Sigma_k$ converges locally smoothly (possibly with multiplicity) to a smooth, almost-embedded $h$-PMC hypersurface $\Sigma_\infty$ in $U$. 
\end{lemma}

\begin{proof}
The proof proceeds via an essentially standard point-picking argument: We first show that the compactness holds on a neighbourhood of any limit point $p=\lim p_k$, where $p_k\in \Sigma_k$. Since $(U, g_k)\to (U, g)$ smoothly, there exists $\rho>0$ so that $B_\rho(p)$ and any sufficiently large $B_\rho(p_k)$ may be identified with a subset $U' \subset \mathbb{R}^{n+1}$ (with their respective metrics), so that $p$ corresponds to $0\in U'$. Moreover, up to taking a subsequence we may assume that $T_{p_k}\Sigma_k \to P$ for some hyperplane $P$, and we may choose the chart so that $P \simeq \{x_{n+1}=0\}$. 

By almost embeddedness and the uniform curvature bounds, there exists $r>0$ (independent of $k$) so that $\Sigma_k\cap B_r(p)$ may be written as a union of ordered graphs \[u_k^{(1)}\leq u_k^{(2)}\leq\cdots \leq u_k^{(m)}\] over $P$. By the uniform area bounds (and the monotonicity formula for bounded first variation), the number of graphs is uniformly bounded. The $u_k^{(i)}$ all have a uniform $C^2$ bound, and satisfy the graphical PMC equation $\mathcal{H}u_k^{(i)} = \sigma_{k,i} h(x,u_k^{(i)}(x)),$ where $\sigma_{k,i}=\pm1$ depending on the orientation of the corresponding sheet. Passing to a subsequence, we may assume that the number of sheets and the signs $\sigma_{k,i}$ are independent of $k$. 

Ascoli's theorem and a diagonal argument gives a subsequence so that for each $i$, $u_k^{(i)}$ converges in $C^{1,\alpha}$ on compact subsets. Schauder theory (and the fact that $h_k\to h$ smoothly) gives that the convergence is in fact smooth. This establishes the compactness on a neighbourhood of any $p=\lim p_k$. 

But the limit set is closed, so covering by the neighbourhoods above, we may use the usual diagonal argument to ensure convergence on any compact subset of $U$. (Note that on any transition regions, ordered sheet decompositions of each $\Sigma_k$ must have consistent orientations, so the orientations of the sheets of $\Sigma_\infty$ are also consistent by the smooth convergence.) 
\end{proof}

We now obtain the compactness in low dimensions for PMC hypersurfaces which are stable under ambient deformations.

\begin{proposition}
Let $2\leq n\leq 6$, and $U\subset M$ be an open subset. If $\Sigma\subset U$ is a smooth, almost embedded, two-sided, stable (under ambient deformations) $h$-PMC hypersurface in $U$ with $\pr \Sigma\cap U=\emptyset$ and $\area(\Sigma)\leq C$, then there exists $C_1$ depending only on $n$, $M$, $c=\sup|h|$, $C$, such that 
\[ |A^\Sigma|^2(x) \leq \frac{C_1}{\dist^2_M(x,\pr U)}, \text{\qquad for all $x\in\Sigma$.}\] 

In particular, if $\Sigma_k\subset U$ is a sequence of smooth, almost embedded, two-sided, stable $h_k$-hypersurfaces in $U$ with $\pr \Sigma_k\cap U=\emptyset$ and $\sup_k \area(\Sigma_k)<\infty,$ and $h_k$ converges locally smoothly to $h$, then up to a subsequence, $\Sigma_k$ converges locally smoothly (possibly with multiplicity) to some $h$-PMC hypersurface $\Sigma$ which is stable (under ambient deformations) in $U$. 
\end{proposition}
\begin{proof}
To prove the curvature estimate, we employ the standard point-picking setup. Recall that $M$ is isometrically embedded in some $\mathbb{R}^N$. 

Suppose there is a sequence of counterexamples, namely $\Sigma_k$ satisfying the hypotheses but with $p_k \in \Sigma_k$ such that \[\rho_k:= \sup |A^{\Sigma_k}| \dist_M(\cdot,U)= |A^{\Sigma_k}|(p_k) \dist_M(p_k,\pr U) \to \infty.\] We may assume $p_k\to p$ in $M$. 

Let $\lambda_k = |A^{\Sigma_k}|(p_k)\to\infty$ and consider the rescaled hypersurfaces $\Sigma_k' = \lambda_k(\Sigma_k - p_k) + p$ in $U'_k = \lambda_k(U_k - p_k) + p$. Note that the $U'_k$ converge smoothly to $T_pM \simeq \mathbb{R}^{n+1}$ on compact subsets. Each $\Sigma'_k$ is a $h'_k:=\lambda_k^{-1}h$ hypersurface in $U'_k$ with $|A^{\Sigma'_k}|(0)=1$, and for sufficiently large $k$, we will have  $\sup_{\Sigma_k'\cap B_{\rho_k/2}} |A^{\Sigma'_k}| \leq 3$. 

Since $\rho_k\to\infty$ and $\lambda_k^{-1}c\to 0$, by Lemma \ref{lem:bdd-cpctness}, the $\Sigma'_k \cap B_{\rho_k/2}$ converge smoothly on compact subsets to a complete, two-sided minimal hypersurface $\Sigma_\infty$ in $T_pM\simeq \mathbb{R}^{n+1}$ (possibly with multiplicity $m$; note that, by the strong maximum principle, an almost embedded \textit{minimal} hypersurface is necessarily embedded). In particular, $|A^{\Sigma_\infty}|(0)=1$. 

We claim that $\Sigma_\infty$ is stable. Indeed, consider $\phi \in C^\infty_c(\Sigma_\infty)$. Then we may extend\footnote{For instance, by the Kirszbraun theorem.} $\phi\nu$ to a Lipschitz vector field $V$ with compact support in a neighbourhood $\mathcal{U}$ of $\supp\phi$ in $\mathbb{R}^N$. 

Let $V_k = \pi_k(V)$ be the tangent projection of $V$ to $TU'_k$, and $\phi_k= \langle V, \nu_{\Sigma'_k}\rangle$. We claim that \begin{equation}\label{eq:lim}\II_{\Sigma_\infty}(\phi, \phi) =m\lim_{k\to\infty} \II_{\Sigma'_k}(\phi_k,\phi_k) .\end{equation}
Indeed, for large $k$, we have $\mathcal{U} \subset B_{\rho_k/2}$, so $U'_k \to T_p M$ and $\Sigma'_k \to \Sigma_\infty$ smoothly in $\mathcal{U}$, and certainly $h_k=\lambda_k^{-1}h \to0$ smoothly. Note that under the local graphical convergence $\Sigma'_k \to \Sigma_\infty$, the projections $\phi_k$ may differ on each sheet of $\Sigma'_k$, but nevertheless the restriction to each sheet must converge smoothly to $\phi$. As in Remark \ref{rmk:2var-cts}, the second variation on each sheet will converge to $II_{\Sigma_\infty}(\phi,\phi)$, which establishes (\ref{eq:lim}) (strictly, one may argue using a partition of unity adapted to a cover of $\supp V$). Using that $\Sigma'_k$ is stable under the ambient deformation by $V_k$, we have that \[ \II_{\Sigma'_k}(\phi_k,\phi_k) \geq 0\] for all $k$, hence $\II_{\Sigma_\infty}(\phi, \phi) \geq 0$ and $\Sigma_\infty$ is stable.

Thus $\Sigma_\infty$ is a complete, two-sided minimal hypersurface in $\mathbb{R}^{n+1}$, and the area bounds imply that $\Sigma_\infty$ has Euclidean area growth. By the stable Bernstein theorem for minimal hypersurfaces, $\Sigma_\infty$ must be flat; but we had $|A^{\Sigma_\infty}|(0)=1$, which gives the desired contradiction. 

With curvature estimates in hand, the final compactness statement follows immediately from Lemma \ref{lem:bdd-cpctness}. (The stability follows similarly to the above: Let $\Sigma$ be the smooth limit and consider an ambient vector field $X$. Stability of $\Sigma_k$ gives $\II_{\Sigma_k}(\phi_k,\phi_k)\geq0$ where $\phi_k = \langle X, \nu_{\Sigma_k}\rangle$, so passing to the limit we have $\II_\Sigma(\phi,\phi)\geq 0$, where $\phi = \langle X, \nu_\Sigma\rangle$ as desired.)
\end{proof}

\section{Good set}
\label{sec:good-set}

In this section, $M$ will denote a closed \Rm\ \mf\ and $h \in L^1(M)$ unless mentioned otherwise. The main result of this section is Proposition \ref{p.good.set.int}, which asserts that a $h$-good set (defined in Definition \ref{d.h.good}) satisfies property $(P_\et)$ (defined in Definition \ref{d.prop.P.eta}). Its proof uses Propositions \ref{p.nested.1} -- \ref{p.nested.3}; their proofs are deferred from the initial presentation to Section \ref{sec:good-set-proofs}.

\begin{defn}\label{d.h.good}
Let $M$ be a complete \Rm\ \mf,  $\Om\subset M$ be an open set with smooth boundary, $h \in L^1(M)$. $\Om$ is said to be a $h$-good set if 
\begin{equation}\label{e.h.good}
4P(\Om)+3V^{|h|}(M\setminus \Om)< W^h(\Om).
\end{equation}	
\end{defn}

\begin{defn}
    A map $\Ph : [0,1] \to \cm$ is called nested if $\Ph(t_1) \subset \Ph(t_2)$ whenever $ t_1 < t_2$.
\end{defn}

\begin{defn}\label{d.minimizer.A^h}
    Let $\ve > 0$, $F,\;F' \in \cm$, $F \subset F'$. $E \in \cm$ is said to be an $(\ve, F, F')$-minimizer for $\ah$ if 
    $$ \ah(E) \le \inf \{\ah(E'): E' \in \cm,\; F \subset E' \subset F'\} + \ve.$$
\end{defn}

\begin{pro}\label{p.nested.1}
	Let $\Ph:[0,1]\to \cm$ be a $\cF$-\cts\ nested map, $E$ be an $(\ve, F, \Ph(1))$-minimizer (where $F \in \cm$ \w\ $ F \subset \Ph(1)$) for $\ah$. Then \tes\ a $\cF$-\cts\ nested map $\tilde{\Ph}:[0,1]\to \cm$ satisfying $\sup_{t\in [0,1]} \ah(\tilde{\Ph}(t)) \leq \sup_{t\in [0,1]} \ah(\Ph(t))+ \ve$, $\tilde{\Ph}(0) \subset \Ph(0)$ and $\tilde{\Ph}(1)=E$. Moreover, if $\Ph(0)\subset E$, one can choose $\tilde{\Ph}(0) = \Ph(0)$.
	
	Similarly, let $\Ph:[0,1]\to \cm$ be a $\cF$-\cts\ nested map, $E$ be an $(\ve, \Ph(0), F)$-minimizer (where $F \in \cm$ \w\ $ \Ph(0) \subset F$) for $\ah$. Then \tes\ a $\cF$-\cts\ nested map $\tilde{\Ph}:[0,1]\to \cm$ satisfying $\sup_{t\in [0,1]} \ah(\tilde{\Ph}(t)) \leq \sup_{t\in [0,1]} \ah(\Ph(t))+ \ve$, $\tilde{\Ph}(0) =E$ and $\Ph(1)\subset \tilde{\Ph}(1)$. Moreover, if $E\subset \Ph(1)$, one can choose $\tilde{\Ph}(1) = \Ph(1)$.
\end{pro}

\begin{pro}\label{p.nested.2}
Let $\Ph:[0,1]\to \cm$ be a $\bF$-\cts\ map. Then \fa\ $\ve>0$, \tes\ a $\cF$-\cts\ nested map $\tilde{\Ph}:[0,1]\to \cm$ \st\ $\sup_{t\in [0,1]} \ah(\tilde{\Ph}(t)) \leq \sup_{t\in [0,1]} \ah(\Ph(t))+\ve$, $\tilde{\Ph}(0)\subset \Ph(0)$ and $\Ph(1)\subset \tilde{\Ph}(1)$.  
\end{pro}

\begin{pro}\label{p.nested.3}
	Let $\Xi \subset M$ be an open set \w\ smooth boundary. For all $\ve>0$, \tes\ $\tilde{\ve}>0$ \st\ the following two conditions are satisfied.
	\begin{enumerate}[(i)]
		\item If $E\in \cm$ satisfies $V(\Xi \cap E)\leq \tilde{\ve}$, then \tes\ $\cF$-\cts\ map $\Ph:[0,1]\to \cm$ \st\ $\Ph(0)=E$, $\Ph(1)\cap \Xi=\emptyset$ and $\sup_{t\in [0,1]} \ah(\Ph(t)) \leq \ah(E)+P(\Xi)+\ve$.
		
		\item If $E\in \cm$ satisfies $V(\Xi\setminus E)\leq \tilde{\ve}$, then \tes\ $\cF$-\cts\ map $\Ph:[0,1]\to \cm$ \st\ $\Ph(0)=E$, $\Ph(1)\supset \Xi$ and $\sup_{t\in [0,1]} \ah(\Ph(t)) \leq \ah(E)+P(\Xi)+\ve$.
	\end{enumerate}	
\end{pro}

The next Lemma \ref{l.isoper} follows from the properties of the isoperimetric profile of the compact manifold $\overline{\Xi}$.
\begin{lem}\label{l.isoper}
Let $\Xi \subset M$ be an open set \w\ smooth boundary. For all $\ve_1>0$, \tes\ $\ve_2>0$ \st\ 
\begin{equation}\label{e.isoper}
\text{if } E\in \cm \text{ satisfies } \min\{V( \Xi \cap E),V(\Xi\setminus E)\}>\ve_1, \text{ then } P(E,\Xi)>\ve_2.
\end{equation}
\end{lem}

\begin{lem}\label{l.cacc.1}
	Let $E,F\in \cm$ and $U \subset M$ be an open set. Then
	\begin{enumerate}[(i)]
		\item $P(E\cap F,U)+P(E \cup F, U)\leq P(E,U)+P(F,U)$;
		\item $|P(E^c\cap F, U)-P(E\cap F , U)|\leq P(F,U)$.
	\end{enumerate}
\end{lem}
\begin{proof}
For the proof of (i), we refer to \cite[Proof of Proposition 2.1]{d}. To prove (ii), we note that for all compactly supported vector-fields $X$ on $U$ \w\ $\|X\|_{C^0}\leq 1$,
\[\int_{E^c\cap F}\operatorname{div}X d\cH^{n+1}=\int_{E\cap F}\operatorname{div}(-X) d\cH^{n+1}+\int_{F}\operatorname{div}X d\cH^{n+1}\le P(E\cap F, U)+P(F,U).\]
\hn\ $P(E^c\cap F, U)\le P(E\cap F, U)+P(F,U)$, which implies (ii).
\end{proof}

\begin{lem}\label{l.cacc.2}
 Let $E,F\in \cm$, $\Xi \subset M$ be an open set with smooth boundary. Then
	$$P\left((E\cap \Xi)\cup (F\cap \Xi^c)\right) \leq P(E,\Xi)+P(F,\overline{\Xi}^c)+P(\Xi).$$
\end{lem}
\begin{proof}\newcommand{\dv}{\operatorname{div}}
Let $\ch:\bbr\ra [0,1]$ be a smooth function \st\ 
\begin{equation}
\ch(t)=\begin{cases}
 0 & \text{ if } t\in (-\infty,1], \\
1 & \text{ if } t\in [2,\infty),
\end{cases}
\end{equation}
and $ \ch'\ge 0$. For $\ve>0$, let $\ch_{\ve}(t)=\ch(t/\ve)$. Let $d:M\ra \bbr$ be a smooth \fn\ \st\ $d$ is equal to the signed distance from $\del \Xi$ on a tubular neighborhood of $\del \Xi$ and $d>0$ on $\Xi^c$. We define $\tht_{\ve}:M\ra [0,1]$ by $\tht_{\ve}(x)=\ch_{\ve}(|d(x)|)$. Then $\tht_{\ve}=\tht_{\ve}^{+}+\tht_{\ve}^{-}$, where $\tht_{\ve}^{+}(x)=\ch_{\ve}(d(x))$, $\tht_{\ve}^{-}(x)=\ch_{\ve}(-d(x))$, $\text{spt}(\tht_{\ve}^+) \Subset \overline{\Xi}^c$, $\text{spt}(\tht_{\ve}^-) \Subset \Xi$. Let $X$ be an arbitrary vector-field on $M$ \w\ $\|X\|_{C^0}\leq 1$ and $R=(E\cap \Xi)\cup (F\cap \Xi^c)$. Then
\begin{align}
\int_R \dv Xd\cH^{n+1}& =\int_R \dv (\tht_{\ve}^-X)d\cH^{n+1}+\int_R \dv (\tht_{\ve}^+X)d\cH^{n+1}+\int_R \dv ((1-\tht_{\ve})X)d\cH^{n+1}\\
& \leq P(E,\Xi) + P(F,\overline{\Xi}^c)+\int_R \dv ((1-\tht_{\ve})X)d\cH^{n+1}.
\end{align}
\tf\ it is enough to show that 
\begin{equation}
\liminf_{\ve\to 0} \int_R \dv ((1-\tht_{\ve})X)d\cH^{n+1} \leq P(\Xi).
\end{equation}
We note that 
$$\dv((1-\tht_{\ve})X)=(1-\tht_{\ve})\dv X + \ch_{\ve}'(-d)\langle \na d, X\rangle - \ch_{\ve}'(d)\langle \na d, X\rangle;$$
using the coarea formula, we obtain 
\begin{multline}
\int_R \dv ((1-\tht_{\ve})X)d\cH^{n+1} = \int_R (1-\tht_{\ve})\dv X d\cH^{n+1}\\
+\int_{1}^{2}\ch'(s)\left(\int_{\{d=-s\ve\}\cap E}\vp d\cH^{n}-\int_{\{d=s\ve\}\cap F}\vp d\cH^{n}\right)ds,
\end{multline}
where $\vp=\langle \na d, X\rangle$. By the dominated convergence theorem,
$$\lim_{\ve\to 0}\int_R (1-\tht_{\ve})\dv X d\cH^{n+1}=0;$$
\hn\ we need to show that
\begin{equation}\label{e.int.ch'}
\liminf_{\ve\to 0}\int_{1}^{2}\ch'(s)\left(\int_{\{d=-s\ve\}\cap E}\vp d\cH^{n}-\int_{\{d=s\ve\}\cap F}\vp d\cH^{n}\right)ds\le P(\Xi).
\end{equation}
Let 
$$\xi_{\ve}(s)=\int_{\{d=-s\ve\}\cap E}\vp d\cH^{n}-\int_{\{d=s\ve\}\cap F}\vp d\cH^{n}.$$
Using the Fermi coordinates of $\del\Xi$, one can show that \te\ $\ta>0$ and a smooth \fn\ $G:\del\Xi\times (-\ta,\ta)\ra (0,\infty)$ \st\ $G(y,0)=1$ and if $\ve$ is sufficiently small,
\[\xi_{\ve}(s)=\int_{\del \Xi}\left(\vp(y,-\ve s)G(y,-\ve s)\ch_E(y,-\ve s)-\vp(y,\ve s)G(y,\ve s)\ch_F(y,\ve s)\right)d\cH^n(y). \]
We note that 
\begin{multline}	
\vp(y,-\ve s)G(y,-\ve s)\ch_E(y,-\ve s)-\vp(y,\ve s)G(y,\ve s)\ch_F(y,\ve s) \\
= \vp(y,0)G(y,0)(\ch_E(y,-\ve s)-\ch_F(y,\ve s))
+(\vp(y,-\ve s)G(y,-\ve s)-\vp(y,0)G(y,0))\ch_E(y,-\ve s)\\
+(\vp(y,0)G(y,0)-\vp(y,\ve s)G(y,\ve s))\ch_F(y,\ve s).
\end{multline}
\hn\
\begin{multline}
\xi_{\ve}(s)\leq P(\Xi)+\int_{\del \Xi}|\vp(y,-\ve s)G(y,-\ve s)-\vp(y,0)G(y,0)|d\cH^n(y)\\
+\int_{\del \Xi} |\vp(y,0)G(y,0)-\vp(y,\ve s)G(y,\ve s)|d\cH^n(y),
\end{multline}
which implies \eqref{e.int.ch'}.
\end{proof}

\begin{lem}\label{l.isoper.2}
Let $\Xi \subset M$ be an open set with smooth boundary, $\de_1,\de_2>0$ \st\ $2\de_1<V(\Xi)$ and \eqref{e.isoper} holds for $(\ve_1,\ve_2)=(\de_1,\de_2)$. Then \fa\ $\Ph\in \widetilde{\sF}(\Xi)$, \te\ $0<a<b<1$ \st\ the following conditions are satisfied.
\begin{enumerate}[(i)]
	\item $P(\Ph(a),\Xi), P(\Ph(b),\Xi)\leq \de_2$;
	\item $V(\Xi\cap \Ph(a)), V(\Xi \setminus \Ph(b))\leq \de_1$;
	\item $\|\del \Ph(t)\|(\overline{\Xi})\geq \de_2$ \fa\ $t\in [a,b]$.
\end{enumerate}
\end{lem}

\begin{proof}
Let \begin{align*}
& I_1=\left\{t\in [0,1]:P(\Ph(t),\Xi)\leq \de_2,\; V(\Xi \cap \Ph(t))\leq \de_1\right\};\\
& I_2=\left\{t\in [0,1]:P(\Ph(t),\Xi)\leq \de_2,\; V(\Xi \setminus \Ph(t))\leq \de_1\right\}.
\end{align*}
Then $I_1, I_2$ are closed subsets of $[0,1]$, $0\in I_1$, $1\in I_2$, $I_1\cap I_2=\emptyset$ and 
\begin{equation}\label{e.I_1.I_2}
I_1 \cup I_2 = \left\{t\in [0,1]:P(\Ph(t),\Xi)\leq \de_2\right\}.
\end{equation}
Let $a=\max I_1$ and $b=\min (I_2\cap [a,1])$. Then $a<b$ and items (i), (ii) of the lemma are satisfied. Moreover, for all $t\in (a,b)$, $P(\Ph(t),\Xi)>\de_2$, which implies $\|\del \Ph(t)\|(\overline{\Xi})\geq \de_2$ \fa\ $t\in [a,b]$.
\end{proof}

 \begin{pro}\label{p.good.set.int}
 Let $M$ be a closed Riemannian manifold, $h \in L^1(M)$ and $\Om\subset M$ be a $h$-good set. Then there exist $\Tht \Subset \Om$, an open set \w\ smooth boundary, and $\et>0$ \st\ the following holds. For all $\Ph\in \widetilde{\sF}(\Om)$, there exists $t_0\in [0,1]$ \st\ $\ah(\Ph(t_0))\geq \widetilde{W}^h(\Om)$ and $\|\del \Ph(t_0)\|(\overline{\Tht})\geq \et$.
 \end{pro}

\begin{proof}\newcommand{\tpts}{\tilde{\Ps}(t_*)}
For $\ta>0$, let
\[\la(\ta)=\sup\left\{V^{|h|}(E):E\subset \Om \text{ is } \cH^{n+1}\text{-measurable and }V(E)\leq \ta\right\}.\]
Then \(\lim\limits_{\ta\to 0}\la(\ta)=0\). Let $\ta_0=\frac{1}{10}P(\Om)$. We choose $0<\ta_1<\frac{1}{2}V(\Om)$ so that $\la(3\ta_1)\leq \ta_0$ and items (i), (ii) of Proposition \ref{p.nested.3} are satisfied for $\Xi= \Om$ and $(\ve,\tilde{\ve})=(\ta_0,3\ta_1)$. Let $\Tht\Subset\Om$ be an open set \w\ smooth boundary \st\ $V(\Tht)> 2\ta_1$, $V(\Om \setminus \Tht)\leq \ta_1$ and $|P(\Om)-P(\Tht)|\leq \ta_0$. Let $0<\ta_2\le \ta_0$ so that \eqref{e.isoper} is satisfied for $\Xi= \Tht$ and $(\ve_1,\ve_2)=(\ta_1,\ta_2)$. We define $\et = \ta_2$.
	
Let us assume by contradiction that there exists $\Ps\in \widetilde{\sF}(\Om)$ for which the assertion does not hold, i.e. for all $t\in [0,1]$ for which $\ah(\Ps(t))\geq \widetilde{W}^h(\Om)$, there holds $\|\del \Ps(t)\|(\overline{\Tht})< \ta_2$. By Lemma \ref{l.isoper.2}, there exist $0\le a<b \le 1$ \st\ 
\begin{align}
&P(\Ps(a),\Tht)=P(\Ps(b),\Tht)\leq \ta_2,\; V(\Tht \cap \Ps(a)), V(\Tht\setminus \Ps(b))\leq \ta_1,\label{e.Ps.a,b}\\
&\|\del \Ps(t)\|(\overline{\Tht})\geq \ta_2\;\forall\;t\in [a,b]\label{e.P.Ps(t)}.
\end{align} 
\eqref{e.P.Ps(t)} implies that there exists $\de>0$ \st
\begin{equation}\label{e.A^h.Ps}
\sup_{t\in [a,b]}\ah(\Ps(t))\leq \widetilde{W}^h(\Om)-2\de.
\end{equation}
By Proposition \ref{p.nested.2}, there exists a nested sweepout $\tilde{\Ps}:[0,1]\ra \cm$ \st\ $\tilde{\Ps}(0)\subset \Ps(a),\; \Ps(b)\subset \tilde{\Ps}(1)$ and
\begin{equation}\label{e.A^h.tldPs}
\sup_{t\in [0,1]}\ah(\tilde{\Ps}(t)) \leq \widetilde{W}^h(\Om)-\de.
\end{equation}

Let $L=\Ps(a)\cup\Ps(b)^c$ so that $L^c=\Ps(a)^c\cap \Ps(b)$ and $\Ps'(t)=\tilde{\Ps}(t)\cup L \cup \Tht^c$, $t\in[0,1]$. Then $\Ps'(0)=L\cup \Tht^c$, $\Ps'(0)\cap \Tht = L\cap \Tht$ and $\Ps'(1)=M$. Hence, $V(\Ps'(0)\cap \Tht)\leq 2\ta_1$, which implies $V(\Ps'(0)\cap \Om)\leq 3\ta_1$. Moreover, by Lemma \ref{l.cacc.2}, Lemma \ref{l.cacc.1} (i) and \eqref{e.Ps.a,b},
\begin{align}
& P(\Ps'(0))\leq P(L,\Tht)+P(\Tht)\leq 2\ta_2+P(\Tht)\le 2\ta_2+\ta_0+P(\Om);\\
& V^{-h}(\Ps'(0))\leq V^{|h|}(\Om^c)+V^{|h|}((L\cap \Tht)\cup (\Om \setminus \Tht))\le V^{|h|}(\Om^c)+\la(3\ta_1).
\end{align}
\tf\ Proposition \ref{p.nested.3} implies that there exists $\Upsilon:[0,1]\ra \cm$ \st\ $\Upsilon(0)=\Ps'(0)$, $\Upsilon(1)\cap \Om = \emptyset$ and 
\begin{equation}\label{e.A^h.Upsilon}
\sup_{t\in [0,1]}\ah(\Upsilon(t))\le \ah(\Ps'(0))+P(\Om)+\ta_0\le 2P(\Om)+V^{|h|}(\Om^c)+2\ta_0+2\ta_2+\la(3\ta_1).
\end{equation}
Concatenating $\Upsilon$ and $\Ps'$ we obtain $\Ps''\in \widetilde{\sF}(\Om)^*$, which implies $\sup_{t \in [0,1]} \ah(\Ps''(t)) \ge \widetilde{W}^h(\Om)$. Since $\Om$ is $h$-good, \eqref{e.W.tldW}, \eqref{e.h.good}, \eqref{e.A^h.Upsilon} imply $$\sup_{t \in [0,1]} \ah(\Upsilon(t)) < \widetilde{W}^h(\Om).$$
So \tes\ $t_*\in [0,1]$ \st\ 
\begin{equation}
\ah(\Ps'(t_*)) > \widetilde{W}^h(\Om) - \de \ge \ah(\tilde{\Ps}(t_*)) \; (\text{by } \eqref{e.A^h.tldPs}).
\end{equation}
Hence, using Lemma \ref{l.cacc.2}, we obtain 
\begin{equation}
P(\Tht)+P(\tpts\cup L,\Tht)-V^h(\tpts\cup L\cup \Tht^c)>P(\tpts,\Tht)+P(\tpts, \overline{\Tht}^c)-V^h(\tpts),
\end{equation}
which implies
\begin{equation}
P(\tpts,\Tht)-P(\tpts\cup L,\Tht)+P(\tpts, \overline{\Tht}^c)<P(\Tht)-V^h(\tpts\cup L\cup \Tht^c)+V^h(\tpts).
\end{equation}
\tf\ by Lemma \ref{l.cacc.1} (i) and \eqref{e.Ps.a,b},
\begin{align}\label{e.P.tpts}
P(\tpts\cap L, \Tht)+P(\tpts,\overline{\Tht}^c) & <P(L,\Tht)+P(\Tht)+V^{|h|}(L\cup \Tht^c)\\
& \le P(\Om)+V^{|h|}(\Om^c)+\ta_0+2\ta_2+\la(3\ta_1).
\end{align}

Let
\begin{align*}
& E^*_0=\left((\tpts\cap L)\cap \Tht \right) \cup (\tpts \cap \Tht^c),\\
& E^*_1 = \left((\tpts\cup L^c)\cap \Tht \right) \cup (\tpts \cap \Tht^c).
\end{align*}
Then $\tilde{\Ps}(0)\subset E^*_0 \subset \tpts \subset E^*_1 \subset \tilde{\Ps}(1)$. Moreover, Lemma \ref{l.cacc.2}, \eqref{e.P.tpts} and \eqref{e.Ps.a,b} imply
\begin{align}\label{e.A^h.E0*}
\ah(E^*_0)& \le P(\tpts\cap L, \Tht)+P(\tpts,\overline{\Tht}^c)+P(\Tht)+V^{|h|}(L\cup \Tht^c)\\
& <2P(\Om)+2V^{|h|}(\Om^c)+2\ta_0+2\ta_2+2\la(3\ta_1).
\end{align}
We assumed that $V^h(M)\ge 0$. \tf\ Lemmas \ref{l.cacc.1}, \ref{l.cacc.2}, \eqref{e.P.tpts} and \eqref{e.Ps.a,b} imply
\begin{align}\label{e.A^h.E1*}
\ah(E^*_1)& \le P(\tpts^c\cap L,\Tht)+ P(\tpts,\overline{\Tht}^c)+P(\Tht) -V^h(E^*_1)+V^h(M)\\
& \le P(\tpts\cap L,\Tht)+P(L,\Tht)+P(\tpts,\overline{\Tht}^c)+P(\Tht)+V^{|h|}(L\cup \Tht^c)\\
& <2P(\Om)+2V^{|h|}(\Om^c)+2\ta_0+4\ta_2+2\la(3\ta_1).
\end{align}

Let $\epsilon = \min \{\de/4, \ta_0\}$, $E^{\bullet}_0,E^{\bullet}_1\in \cm$ \st\ $E^{\bullet}_0$ is an $(\epsilon, \tilde{\Ps}(0), E^*_0)$-minimizer for $\ah$ and $E^{\bullet}_1$ is an $(\epsilon, E^*_1, \tilde{\Ps}(1))$-minimizer for $\ah$. Using Proposition \ref{p.nested.1} twice and \eqref{e.A^h.tldPs}, we construct a nested sweepout $\hat{\Ps}:[0,1]\ra \cm$ \st\ $\hat{\Ps}(0)=E^{\bullet}_0$, $\hat{\Ps}(1)=E^{\bullet}_1$ and
\begin{equation}\label{e.A^h.Ps.hat}
\sup_{t\in [0,1]}\ah(\hat{\Ps}(t))\leq \widetilde{W}^h(\Om)- \frac{\de}{2}.
\end{equation}
\eqref{e.Ps.a,b} implies
\begin{equation}
V(\Om\cap E^{\bullet}_0 )\leq V(\Om \cap E^*_0 )\le 3\ta_1,\; V(\Om\setminus E^{\bullet}_1 )\leq V(\Om \setminus E^*_1 )\le 3\ta_1.
\end{equation}
Hence, by Proposition \ref{p.nested.3}, there exists $\Upsilon_i:[0,1]\ra \cm$, $i=0,1$, \st\ $\Upsilon_i(0)=E^{\bullet}_i$, $\Om \cap \Upsilon_0(1) = \emptyset$, $\Om\setminus \Upsilon_1(1) = \emptyset$ and (by \eqref{e.A^h.E0*}, \eqref{e.A^h.E1*})
\begin{equation}\label{e.A^h.Up_i}
\sup_{t\in [0,1]}\ah(\Upsilon_i(t))\le \ah(E^{*}_i)+P(\Om)+ \ta_0<4P(\Om)+2V^{|h|}(\Om^c).
\end{equation}
Concatenating $\Upsilon_1,\;\Upsilon_2$ and $\hat{\Ps}$, we obtain $\Upsilon' \in \widetilde{\sF}(\Om)^*$. Since $\Om$ is $h$-good, \eqref{e.A^h.Ps.hat}, \eqref{e.A^h.Up_i}, \eqref{e.h.good} and \eqref{e.W.tldW} imply that
\[\sup_{t\in [0,1]}\ah(\Upsilon'(t)) < \widetilde{W}^h(\Om),\]
which contradicts \eqref{e.tilde.W.*}.
\end{proof}

%%%%%%%%%%%%%%%%%%%%%%%%%%%%%%%%%%%%%%%%%%%%%%%%%%%%%%%%%%%%
%%%%%%%%%%%%%%%%%%%%%%%%%%%%%%%%%%%%%%%%%%%%%%%%%%%%%%%%%%%%
\subsection{Proofs of Propositions}
%%%%%%%%%%%%%%%%%%%%%%%%%%%%%%%%%%%%%%%%%%%%%%%%%%%%%%%%%%%%
%%%%%%%%%%%%%%%%%%%%%%%%%%%%%%%%%%%%%%%%%%%%%%%%%%%%%%%%%%%%
\label{sec:good-set-proofs}

In this subsection, we detail the proofs of Propositions \ref{p.nested.1}, \ref{p.nested.2}, and \ref{p.nested.3}.

%%%%%%%%%%%%%%%%%%%%%%%%%%%%%%%%%%%%%%%%%%%%%%%%%%%%%%%%%%%%
\subsubsection{Proof of Proposition \ref{p.nested.1}}
%%%%%%%%%%%%%%%%%%%%%%%%%%%%%%%%%%%%%%%%%%%%%%%%%%%%%%%%%%%%
    This statement follows the procedure outlined in the section ``Gluing two nested families" on page 200 of \cite{cl}.
    One should also refer to Figure 1 on page 188 of the same article.
    
    To prove the first part of the proposition we define  $$\tilde{\Phi}(t) = E \cap \Phi(t)$$
    for each $t$. Clearly every such set is a Caccioppoli set, $\tilde{\Phi}$ is $\mathcal{F}$-continuous, and $\tilde{\Phi}$ satisfies the desired inclusion-exclusion
    properties. To prove the bounds on $\mathcal{A}^h$, we will prove that $\ah(\Ph(t)\cap E)\leq \ah(\Ph(t)) + \ve $, which is equivalent to 
    \begin{equation}\label{e.A^h.E.ph(t).1}
        V^h(\Ph(t)\setminus E)\leq P(\Ph(t))-P(\Ph(t)\cap E)+\ve.
    \end{equation}
    Since $F\subset E\subset \Ph(t)\cup E \subset \Ph(1)$ and $E$ is an $(\ve, F, \Ph(1))$-minimizer for $\ah$, there holds $\ah(\Ph(t)\cup E)\geq \ah(E)-\ve$, which is equivalent to
    \begin{equation}
        \label{e.A^h.E.ph(t).2}
        V^h(\Ph(t)\setminus E)\leq P(\Ph(t)\cup E)-P(E)+\ve.
    \end{equation}
    \eqref{e.A^h.E.ph(t).2} and Lemma \ref{l.cacc.1} (i) imply \eqref{e.A^h.E.ph(t).1}.
    
    To prove the second part of the proposition we define  $$\tilde{\Phi}(t) = E \cup \Phi(t)$$ for each $t$. Clearly every such set is a Caccioppoli set, $\tilde{\Phi}$ is $\mathcal{F}$-continuous, and $\tilde{\Phi}$ satisfies the desired inclusion-exclusion
    properties.  To prove the bounds on $\mathcal{A}^h$, we will prove that $\ah(\Ph(t)\cup E)\leq \ah(\Ph(t)) + \ve $, which is equivalent to 
    \begin{equation}\label{e.A^h.E.ph(t).3}
        V^h(E \setminus \Ph(t))\geq P(\Ph(t)\cup E)-P(\Ph(t))-\ve.
    \end{equation}
    Since $\Phi(0) \subset \Ph(t)\cap E \subset E \subset F$ and $E$ is an $(\ve, \Ph(0), F)$-minimizer for $\ah$, there holds $\ah(\Ph(t)\cap E)\geq \ah(E)-\ve$, which is equivalent to
    \begin{equation}
        \label{e.A^h.E.ph(t).4}
        V^h(E\setminus \Ph(t))\geq P(E) -P(\Ph(t)\cap E) -\ve.
    \end{equation}
    \eqref{e.A^h.E.ph(t).4} and Lemma \ref{l.cacc.1} (i) imply \eqref{e.A^h.E.ph(t).3}.

%%%%%%%%%%%%%%%%%%%%%%%%%%%%%%%%%%%%%%%%%%%%%%%%%%%%%%%%%%%%
\subsubsection{Proof of Proposition \ref{p.nested.2}}
%%%%%%%%%%%%%%%%%%%%%%%%%%%%%%%%%%%%%%%%%%%%%%%%%%%%%%%%%%%%

Given $E\in \cm$, let $\db{E}$ denote the mod 2 current defined by $E$. $\bM$ denotes the mass norm on the space of mod 2 currents. 

\begin{lem}\label{l.nested}
Let $E_0,E_1\in \cm$, $\de>0$ \st\ $E_0\subset E_1$,
\begin{equation}\label{e.mass.diff.1}
V^{|h|}(E_1\setminus E_0)\leq \de,\quad \bM(\del \db{E_0}-\del \db{E_1})\leq \de.
\end{equation}
Let us fix a $\cF$-\cts\ nested map $\Ps:[0,1]\ra \cm$, $\Ps(0)=\emptyset$, $\Ps(1)=M$, with no concentration of mass. Let $R>0$ so that \fa\ $p\in M$, 
\begin{equation}\label{e.no.conc.mass}
P(E_0,B(p,2R))\leq \de,\quad P(E_1,B(p,2R))\leq \de,\quad \sup_{t\in [0,1]}P(\Ps(t),B(p,2R))\leq \de.
\end{equation}
\sps\ $M$ can be covered by $I$ balls of radius $R$,
\begin{equation}\label{e.flat.diff.1}
V(E_1\setminus E_0)\leq \frac{\de R}{I}.
\end{equation}
Then \tes\ a $\cF$-\cts\ nested map $\Ph:[0,1]\ra \cm$ \w\ $\Ph(0)=E_0, \Ph(1)=E_1$ and 
\begin{equation}
\sup_{t\in [0,1]}\ah(\Ph(t))\leq \min \{\ah(E_0),\ah(E_1)\}+11\de.
\end{equation} 
\end{lem}

\begin{proof}
Let 
\begin{equation}\label{e.cover.I.balls}
M=\bigcup_{i=1}^I B(p_i,R).
\end{equation}
We define a sequence \(\{F_i\}_{i=1}^I\), $F_i\in \cm$ \w
\[E_0=F_0\subset F_1\dots\subset F_I=E_1\]
as follows. Set $F_0=E_0$. \sps\ $F_k$ has been defined for $1\leq k \leq i-1$ and 
\[E_0=F_0\subset F_1\dots\subset F_{i-1}\subset E_1.\]
Then $V(E_1\setminus F_{i-1})\leq V(E_1\setminus E_0)\leq \de R/I$. Using the fact that $\cH^n(\del E) < \infty$ for $E \in \cm$, \cite[Lemma 5.2]{eg} and the coarea formula, one can show that \tes\ $r_i\in (R,2R)$ \st\ the following conditions are satisfied. (In the following lines $B_i$ stands for $B(p_i,r_i)$ and $\bar{B}_i$ stands for the closure of $B(p_i,r_i)$.) 
\begin{equation}
\cH^n(\del B_i\cap \del E_0) = \cH^n(\del B_i\cap \del E_1) = \cH^n(\del B_i\cap \del F_{i-1})=0,
\end{equation}
the reduced boundary of $F_{i-1}\cup (E_1\cap B_i)$ is 
\begin{equation}
(\bar{B}_i^c\cap \del F_{i-1})\cup (B_i\cap \del E_1)\cup ((E_1\setminus F_{i-1})\cap \del B_i)
\end{equation}
and
\begin{equation}
\cH^n((E_1\setminus F_{i-1})\cap \del B_i)\leq \frac{\de}{I}.
\end{equation}
Define $F_i=F_{i-1}\cup (E_1\cap B_i)$. Then $F_{i-1}\subset F_i \subset E_1$. Using induction, one can prove that 
\begin{equation}\label{e.F_i.F_i-1}
F_i = \left(E_0\cap \overline{G_i^0}\right) \cup \left(E_1\cap G_i^1\right)\;\text{ and }\; \del F_i = (G_i^0\cap \del E_0)\cup (G_i^1\cap \del E_1) \cup \left(\cup_{k=1}^i S_{k,i}\right),
\end{equation}
where $G_i^0=\left(\cup_{k=1}^i\bar{B}_k\right)^c$, $G_i^1 = \cup_{k=1}^i B_k$, $S_{k,i}\subset (E_1\setminus F_{k-1})\cap \del B_k$. As a consequence, $\cH^n(S_{k,i})\leq \de/I$ and
\begin{align}\label{A^h.F_i}
\ah(F_i) & =P(E_0,G_i^0)-V^h\left(E_0\cap \overline{G_i^0}\right)+P(E_1,G_i^1)-V^h(E_1\cap G_i^1)+\sum_{k=1}^{i}\cH^n(S_{k,i})\\
& \leq \min\{\ah(E_0),\ah(E_1)\}+\bM(\del\db{E_1}-\del\db{E_0})+V^{|h|}(E_1\setminus E_0)+ \frac{\de i}{I}\\
& \leq \min\{\ah(E_0),\ah(E_1)\}+3\de.
\end{align}
Moreover, \fa\ $p\in M$, 
\begin{align}\label{e.P.F_i.B}
P(F_i,B(p,2R))\leq P(E_0, B(p,2R))+P(E_1, B(p,2R))+\sum_{k=1}^{i}\cH^n(S_{k,i})\leq 3\de.
\end{align}

For each $1\leq i\leq I$, define $\Ph_i:[0,1]\ra \cm$ by $\Ph_i(t)=F_{i-1}\cup (F_i\cap \Ps(t))$. Then $\Ph_i$ is a nested, $\cF$-\cts\ map, $\Ph_i(0)=F_{i-1}$, $\Ph_i(1)=F_i$. We note that
\begin{equation}
\Ph_i(t)=F_{i-1}\cup (E_1\cap B_i\cap \Ps(t)).
\end{equation}
\tf\ $\Ph_i(t)\cap B_i^c = F_{i-1}\cap B_i^c$; $\db{\del \Ph_i}\mres \bar{B}_i^c=\db{\del F_{i-1}}\mres \bar{B}_i^c$;
\begin{align}
P(\Ph_i(t)) & = \|\del \Ph_i(t)\|(\bar{B}_i^c) + \|\del \Ph_i(t)\|(\bar{B}_i)\\
& \leq P(F_{i-1})+7\de \text{ (by Lemma \ref{l.cacc.1} (i) and \eqref{e.P.F_i.B})};\\
|V^h(\Ph_i(t)) &-V^h(F_{i-1})|=|V^h(E_1\cap B_i\cap \Ps(t)\cap F_{i-1}^c)|\leq V^{|h|}(E_1\setminus E_0)\leq \de.
\end{align}
Thus
\[\ah(\Ph_i(t))\leq \ah(F_{i-1})+8\de\leq \min\{\ah(E_0),\ah(E_1)\}+11\de.\]
The map $\Ph$ is obtained by concatenating all the maps $\{\Ph_i\}_{i=1}^I$. 
\end{proof}
\newcommand{\Ups}{\Upsilon}
\begin{lemma}\label{l.gluing}
    Let $\de > 0$, \(\Ups_1, \Ups_2: [0,1] \to \cm\) be $\cF$-\cts\ nested maps \st\ $\Ups_2(0) \subset \Ups_1(1)$. Then \tes\ a $\cF$-\cts\ nested map $\Ups: [0,1] \to \cm$ \st\ $\Ups(0) \subset \Ups_1(0)$, $\Ups(1) \supset \Ups_2(1)$ and for $i=1,2$,
    \begin{equation}
        \sup_{t \in [0,1]} \ah(\Ups(t)) \le \sup_{t \in [0,1]} \ah(\Ups_i(t)) + \de.
    \end{equation}
\end{lemma}

\begin{proof}
Let \(E\) be a $(\de, \Ups_2(0), \Ups_1(1))$-minimizer for \(\ah\). Then by Proposition \ref{p.nested.1}, \tes\ a nested map $\tilde{\Ups}_1: [0,1] \to \cm$ \st\ $\tilde{\Ups}_1(0) \subset \Ups_1(0)$, $\tilde{\Ups}_1(1) = E$, and
\begin{equation}
        \sup_{t \in [0,1]} \ah(\tilde{\Ups}_1(t)) \le \sup_{t \in [0,1]} \ah(\Ups_1(t)) + \de.
    \end{equation}
Similarly, \tes\ a nested map $\tilde{\Ups}_2: [0,1] \to \cm$ \st\ $\tilde{\Ups}_2(0) = E$, $\tilde{\Ups}_2(1) \supset \Ups_2(1)$, and
\begin{equation}
        \sup_{t \in [0,1]} \ah(\tilde{\Ups}_2(t)) \le \sup_{t \in [0,1]} \ah(\Ups_2(t)) + \de.
    \end{equation}
The map $\Ups$ is obtained by concatenating $\tilde{\Ups}_1$ and \(\tilde{\Ups}_2\).
\end{proof}

\begin{proof}[Proof of Proposition \ref{p.nested.2}]
Let $\de_0=\ve/13$. By the interpolation results, \tes\ a $\cF$-\cts\ map $\Ph':[0,1]\ra \cm$, \st\ $\Ph'(0)=\Ph(0)$, $\Ph'(1)=\Ph(1)$, $t \mapsto \del\db{\Ph'(t)}$ is \cts\ \wrt\ the $\bM$ norm and
\begin{equation}\label{e.A^h.Ph'}
\sup_{t\in [0,1]} \ah(\Ph'(t))\leq \sup_{t\in [0,1]} \ah(\Ph(t))+\de_0.
\end{equation} 
Let us fix a $\cF$-\cts\ nested map $\Ps_0:[0,1]\ra \cm$, $\Ps_0(0)=\emptyset$, $\Ps_0(1)=M$, with no concentration of mass. (For example, one can construct such a map by considering the sublevel sets of a Morse function on $M$.) \Tes\ $R_0>0$ \st\ \fa\ $p\in M$ and $t\in [0,1]$,
\begin{equation}\label{e.no.conc.mass.1}
P(\Ph'(t),B(p,2R_0))\leq \de_0,\quad P(\Ps_0(t),B(p,2R_0))\leq \de_0.
\end{equation}
\sps\ $M$ can be covered by $I_0$ balls of radius $R_0$. Let $m\in \bbn$ \st\ $|t_1-t_2|\leq 1/m$ implies 
\begin{equation}\label{e.mass.diff.2}
V^{|h|}(\Ph'(t_1)\triangle \Ph'(t_2))\leq \de_0;\quad V(\Ph'(t_1)\triangle \Ph'(t_2))\leq \frac{\de_0R_0}{I_0};\quad \bM(\del \db{\Ph'(t_1)}-\del \db{\Ph'(t_2)})\leq \de_0.
\end{equation}
For $i=0,\dots,m$, let $F_i=\Ph'(i/m)$. 

We define the \seq\  \(\{D_i\}_{i=0}^{2m-1}\), $D_i\in \cm$, by
\[D_{2i}=F_i\cap F_{i+1},\quad  D_{2i+1}=F_i\cup F_{i+1}.\]
Then $D_{2i}\subset D_{2i+1}$, $D_{2i+1}\setminus D_{2i} = F_i\triangle F_{i+1}$,
\begin{equation}
\del(\db{D_{2i}}-\db{D_{2i+1}})=\del \db{F_i\triangle F_{i+1}}=\del(\db{F_i}-\db{F_{i+1}}).
\end{equation}
Moreover, by Lemma \ref{l.cacc.1} (i),
\begin{equation}
 P(D_{2i})+P(D_{2i+1})\leq P(F_i)+ P(F_{i+1}),\quad V^h(D_{2i})+V^h(D_{2i+1})=V^h(F_i)+ V^h(F_{i+1});
\end{equation}
hence
\begin{equation}\label{e.A^h.D_i}
\ah(D_{2i})+\ah(D_{2i+1})\leq \ah(F_i)+ \ah(F_{i+1}).
\end{equation}
\tf\ by \eqref{e.A^h.D_i}, \eqref{e.A^h.Ph'} and Lemma \ref{l.nested}, \tes\  $\cF$-\cts\ nested map $\tilde{\Ph}_i:[0,1]\ra \cm$ \st\ $\tilde{\Ph}_i(0)=D_{2i},\; \tilde{\Ph}_i(1)=D_{2i+1}$ and 
\begin{equation}\label{e.A^h.Tph_i.t}
\sup_{t\in [0,1]}\ah(\tilde{\Ph}_i(t))\leq \sup_{t\in [0,1]} \ah(\Ph(t))+12\de_0.
\end{equation}
We note that $\tilde{\Ph}_{i+1}(0)\subset F_{i+1}\subset \tilde{\Ph}_i(1)$. In particular, 
\begin{equation}\label{e.tPh.0.0}
    \tilde{\Ph}_0(0) \subset F_0 = \Ph(0),\quad \tilde{\Ph}_{m-1}(1) \supset F_m = \Ph(1).
\end{equation}

The nested map $\tilde{\Ph}$ in Proposition \ref{p.nested.2} is obtained by gluing the maps \(\{\tilde{\Ph}_i\}_{i=0}^{m-1}\) (using Lemma \ref{l.gluing}) in the following way. Applying Lemma \ref{l.gluing} to \(\tilde{\Ph}_0, \tilde{\Ph}_1\) and using \eqref{e.A^h.Tph_i.t}, we obtain a nested map $\hat{\Ph}_1$ \w\ $\hat{\Ph}_1(0) \subset \tilde{\Ph}_0(0)$, $\hat{\Ph}_1(1) \supset \tilde{\Ph}_1(1)$ and
\begin{equation}
    \sup_{t\in [0,1]}\ah(\hat{\Ph}_1(t))\leq \sup_{t\in [0,1]} \ah(\Ph(t))+12\de_0 + \frac{\de_0}{m}.
\end{equation}
Inductively, using Lemma \ref{l.gluing}, one can construct a sequence \(\{\hat{\Ph}_i\}_{i=1}^{m-1}\) of nested maps \st\ for $i >1$,  \(\hat{\Ph}_i(0) \subset \hat{\Ph}_{i-1}(0)\), $\hat{\Ph}_i(1) \supset \tilde{\Ph}_i(1)$ and
\begin{equation}\label{e.A^h.hat.Ph_i}
    \sup_{t\in [0,1]}\ah(\hat{\Ph}_i(t))\leq \sup_{t\in [0,1]} \ah(\Ph(t))+12\de_0 + \frac{i\de_0}{m}.
\end{equation}
Finally, we define $\tilde{\Ph} = \hat{\Ph}_{m-1}$. \eqref{e.tPh.0.0} and \eqref{e.A^h.hat.Ph_i} imply that $\tilde{\Ph}$ satisfies the required properties.
\end{proof}
\begin{comment}
To construct the nested map $\tilde{\Ph}$, we use Proposition \ref{p.nested.2}.  Choose an $(\ve,\tilde{\Ph}_0(1),\tilde{\Ph}_1(0))$-minimizer $E$,  we can apply the first part of Proposition \ref{p.nested.2} to the nested family $\tilde{\Ph}_0$ and $E$, and we can apply the second part of the proposition to the nested family $\tilde{\Ph}_1$ and $E$, gluing together the two resulting families to obtain a new nested family $\tilde{\Ph}_{0,1}$.  By Property 9 above and Proposition \ref{p.nested.2},
    $$ \tilde{\Ph}_{0,1}(0) \subset \tilde{\Ph}_0(0) $$
and
    $$ \tilde{\Ph}_1(1) \subset \tilde{\Ph}_{0,1}(1). $$

Thus, $\tilde{\Ph}_2(0) \subset \tilde{\Ph}_{0,1}(1)$, and so we can apply Proposition \ref{p.nested.2}
to $\tilde{\Ph}_2$ and $\tilde{\Ph}_{0,1}$.  We continue this process until we have merged all of the families $\tilde{\Ph}_k$ together.  At each stage, we observe that the initial set of the new family lies within the initial set of the previous family.  Furthermore, at each stage, the final set of the new family contains the final set in the second family involved in the merging.
    
Combining these two facts, the resulting nested family starts within $\tilde{\Ph}_0(0)$, and ends outside of $\tilde{\Ph}_{m-1}(1)$, and has the desired bounds on $\mathcal{A}^h$.  This completes the proof.
\end{comment}

%%%%%%%%%%%%%%%%%%%%%%%%%%%%%%%%%%%%%%%%%%%%%%%%%%%%%%%%%%%%
\subsubsection{Proof of Proposition \ref{p.nested.3}}
%%%%%%%%%%%%%%%%%%%%%%%%%%%%%%%%%%%%%%%%%%%%%%%%%%%%%%%%%%%%

The proof of the first part is very similar to the proof of
the second part; we will prove the first part and then make some comments about how the proofs differ.

The proof of part (i) follows the proof of Lemma~7.1 in \cite{cl}.  To begin, we choose $r_0 > 0$ sufficiently small so that every ball $B$ of radius $r \in (0,r_0]$ centered at a point in $\overline{\Xi}$ is $2$-biLipschitz diffeomorphic to a ball of the same radius in Euclidean space.  Since $\overline{\Xi}$ is compact, this is possible.

Choose a covering $\{B_i\}$ of $\overline{\Xi}$ by $N$ balls of radius $r_0$ so that concentric balls of radius $r_0/2$, denoted by $1/2 B_i$, still cover $\overline{\Xi}$.  Since $\Xi$ is precompact, we can find some $\eta > 0$ so that, for any measurable subset $U$ of $\Xi$ of volume at most $\eta$, $$ \int_U |h| dV < \frac{\ve}{10N}.$$  Now, choose
$$ \tilde{\ve} = \min \left\{ \frac{\ve r_0}{20 N}, \left(\frac{\ve}{2^{n+1}}\right)^{\frac{n+1}{n}}, \eta \right\}.$$
By the coarea formula, we can choose a covering $\{B'_i\}$
of $\Xi$ by $N$ balls of radius $r_i \in (r_0/2,r_0)$ so
that $$\mathcal{H}^n((\partial B'_i) \cap (E \cap \Xi)) \leq \frac{4 \tilde{\ve}}{r_0} \leq \frac{\ve}{5N}.$$

We now inductively push $E$ outside of $B_1' \cap \Xi, \dots, B_N' \cap \Xi$.  To push $E$ outside of $B_1' \cap \Xi$,
we consider the $2$-biLipschitz map $f : B_1' \rightarrow B$,
where $B$ has the same radius as $B_1'$ and is in $\mathbb{R}^{n+1}$.  Then we can contract $f(B_1' \cap \Xi \cap E)$ to
a single point in $\mathbb{R}^n$ linearly; let this
map be $\tilde{\Phi}$.  To build $\Phi$, we take
$f^{-1} \circ \tilde{\Phi}$.  We modify this map so that
it removes $f^{-1}(0)$ from the final set.

This $\Phi$ is clearly $\mathcal{F}$-continuous.  We examine
how $\mathcal{A}^h$ changes as we move along this path.
Since the sets are nested, $V^h$ can only increase by at most $$ \int_{E \cap \Xi \cap B_1'} |h| dV < \frac{\ve}{10N}.$$
Furthermore, the perimeter can increase by at most $\mathcal{H}^n(B_1' \cap \partial E) + \mathcal{H}^n(\partial \Xi \cap B_1' \cap E) + 2^n \tilde{\ve}^{\frac{n}{n+1}}$.  From our choice of $\tilde{\ve}$, we can conclude that this is bounded by
$\mathcal{H}^n(B_1' \cap \partial \Xi) + \frac{\ve}{2 N}$. Thus, $\mathcal{A}^h$ can increase by at most 
$\mathcal{H}^n(B_1' \cap \partial \Xi) + \frac{3\ve}{5N}.$

We repeat this procedure iteratively for $B_2', B_3', \dots, B_N'$.  After the $i$th step, we have that
$\mathcal{A}^h$ can increase by at most 
$$ \mathcal{H}^n \Big( (\cup_{j=1}^i B_i') \cap \partial \Xi \Big) + \frac{3i \ve}{5N}, $$
from which the conclusion is obvious.

The proof of part (ii) is the same, except that instead of shrinking portions of $E \cap \Xi$, we shrink portions of
$\Xi - E$.

%%%%%%%%%%%%%%%%%%%%%%%%%%%%%%%%%%%%%%%%%%%%%%%%%%%%%%%%%%%%%%%%%%%%%%%%%%%%%%%%%%%%%%%%%%%%%%%%%%%%%%%%
\section{Proof of the main theorem}
\label{sec:noncompact}
%%%%%%%%%%%%%%%%%%%%%%%%%%%%%%%%%%%%%%%%%%%%%%%%%%%%%%%%%%%%%%%%%%%%%%%%%%%%%%%%%%%%%%%%%%%%%%%%%%%%%%%%
In this section we prove Theorem \ref{main theorem}.

\begin{definition}\label{d.W^h(M)>0}
Let $M$ be a complete Riemannian manifold. Let $\cS$ denote the set of all sequences $\{\Om_i\}_{i=1}^{\infty}$ satisfying $\Om_i \subset M$ is a bounded open set with smooth boundary, $\Om_i \subset \Om_{i+1}$ and $M=\cup_{i=1}^{\infty}\Om_i$. For $h\in L^{1}(M)$, define the width of $M$ \wrt\ the $\ah$ functional 
$$W^h(M) = \inf_{\{\Om_i\} \in \cS} \liminf_{i\to \infty} W^h(\Om_i).$$
\end{definition}

\begin{defn}\label{d.prop.P}
Let $M$ be a complete Riemannian manifold. $M$ is said to satisfy property $(\mathbf{P})$ if 
\begin{itemize}
    \item[$(\mathbf{P})$] there exists a nested sequence of bounded open sets \(\Om_i \subset M\) with smooth boundary \st\ \(\lim_{i\to \infty}P(\Om_i)=0\) and $M=\cup_{i=1}^{\infty}\Om_i$.
\end{itemize}
\end{defn}

The following lemma follows from Definitions \ref{d.h.good}, \ref{d.W^h(M)>0}, \ref{d.prop.P}.

\begin{lem}\label{l.existence.of.good.sets}
	Let $M$ be a complete Riemannian manifold which satisfies property $(\mathbf{P})$, $h\in L^1(M)$ and $W^h(M)>0$. Then there exists $\Om \subset M$, which is $h$-good.
\end{lem}

\begin{defn}\label{d.width}
Let $M$ be a complete Riemannian manifold. The width of $M$
\[W(M)=\sup\{W(\Xi):\Xi\subset M \text{ is a bounded open set with smooth boundary}\},\]
where $W(\Xi)$ denotes the relative width of $\Xi$ with respect to the area functional.
\end{defn}

\begin{lem}\label{l.conditions.for.positive.W^h.1}
	Let $M$ be a complete Riemannian manifold, $h\in L^{1}(M)$ \st\ $\|h\|_{L^1} = V^{|h|}(M)<W(M)$. Then $W^h(M) \ge W(M)-V^{|h|}(M)$.
\end{lem}

\begin{proof}
Let $ 0 < \de < W(M)-V^{|h|}(M)$. There exists \(\Om_0\subset M\), a bounded open set with smooth boundary, \st\ \(W(\Om_0)\geq W(M)-\de\). Let $\Om$ be a bounded open set with smooth boundary \st\ \(\Om_0\subset \Om\). Then \(W(\Om)\geq W(\Om_0)\). For all $\Ph\in \sF(\Om)$, \tes\ $t_0\in [0,1]$ \st\ $P(\Ph(t_0),\Om)\geq W(\Om)$; hence $\ah(\Ph(t_0),\Om)\geq W(\Om) - V^{|h|}(\Om)$. Therefore, 
$$W^h(\Om)\geq  W(\Om) - V^{|h|}(\Om)\geq W(M)- \de - V^{|h|}(M),$$
which implies $W^h(M) \ge W(M)-V^{|h|}(M) -\de$. Letting $\de \to 0$, we obtained the required inequality.
\end{proof}

\begin{lem}\label{l.conditions.for.positive.W^h.2}
Let $M$ be a complete Riemannian manifold of finite volume, $\sI$ be the isoperimetric profile of $M$, $h\in L^{\infty}(M)$ \st\ 
\[\nm{h}_{L^{\infty}}<[\sI]_0:=\sup_{0<v<V(M)}\frac{\sI(v)}{v}.\]
Then $W^h(M)>0$.
\end{lem}

\begin{proof}
Let $k=\frac{1}{2}(\nm{h}_{L^{\infty}}+[\sI]_0)$. Then \tes\ $v_0\in (0,V(M))$ \st\ $\sI(v_0)>kv_0$. Fix $ 0 <\ve < \min\{v_0 , V(M)-v_0\}$. It was proved in \cite{cl} that $M$ satisfies property $(\mathbf{P})$ in Definition \ref{d.prop.P}. Hence, \tes\ \(\Om_0\subset M\), a bounded open set with smooth boundary, \st\ $V(M\setminus \Om_0)\leq \ve$ and $P(\Om_0)\leq \ve$. Let $\Om$ be a bounded open set with smooth boundary \st\ \(\Om_0\subset \Om\) and $\Ph\in \mathscr{F}(\Om)$. Then \tes\ $t_0 \in [0,1]$ \st\ $V(\Ph(t_0)\cap \Om_0)=v_0$, which implies $P(\Ph(t_0)\cap \Om_0)\geq \mathscr{I}(v_0)>kv_0$. Using Lemma \ref{l.cacc.2}, we obtain
\begin{equation}
	P(\Ph(t_0),\Om)\geq P(\Ph(t_0),\Om_0)\geq P(\Ph(t_0)\cap \Om_0) - P(\Om_0)>kv_0-\ve. 
\end{equation}
\tf\
\begin{align}
	\ah(\Ph(t_0),\Om) & =P(\Ph(t_0),\Om)-V^h(\Ph(t_0),\Om)\\
	& >kv_0-\ve -\|h\|_{L^{\infty}}(v_0+V(M\setminus \Om_0))\\
	& \geq \frac{1}{2}([\sI]_0-\|h\|_{L^{\infty}})v_0-\ve(1+\|h\|_{L^{\infty}}).
\end{align}
\hn\ 
$$W^h(\Om)\geq \frac{1}{2}([\sI]_0-\|h\|_{L^{\infty}})v_0-\ve(1+\|h\|_{L^{\infty}}).$$
This implies (after letting $\ve\to 0$)
\[W^h(M)\geq \frac{1}{2}([\sI]_0-\|h\|_{L^{\infty}})v_0 > 0.\]
\end{proof}

\begin{rmk}
By \cite[Section 5]{Gut}, for every $n\geq 1$, 
$$\sup\{W(N): N^{n+1} \text{ is a closed Riemannian manifold with } \vol(N)=1\}=\infty.$$ Hence, there exist a sequence of closed Riemannian manifolds $N_i$ \st\ $\sum_{i=1}^{\infty} \vol(N_i) <\infty$ and $\lim_{i\to \infty}W(N_i)= \infty$. By considering a connected sum of the $N_i$'s, one may construct a complete Riemannian manifold $M$ such that $\vol(M)<\infty$ and $W(M)=\infty$. 

For a complete Riemannian manifold with finite volume, the isoperimetric profile $\sI: [0,\vol(M)]\ra [0,\infty)$ is a continuous function, $\sI(0)=0$ and $\sI(v)>0$ for all $v\in (0,\vol(M))$. Hence $[\sI]_0$, the Lipschitz seminorm of $\sI$ at $0$, is always positive. For a closed Riemannian manifold, $[\sI]_0$ is always $\infty$. By \cite[Example 11.10]{str}, there exist complete non-compact Riemannian manifolds with finite volume for which $[\sI]_0=\infty$.
\end{rmk}

\begin{rmk}\newcommand{\rhe}{\rho_{\ve}}
Let $\rh: [0, \infty ) \to [0,1]$ be defined by $\rh(t) = \min\{t, e^{-t}\}$. \Tes\ unique $t_0 \in (0, 1)$ \st\ $t_0 = e^{-t_0}$. For $0 < \ve < 1-t_0$, \tes\ a smooth \fn\ $\rhe: [0, \infty ) \to [0,1]$ \st\ $\rhe = \rh$ outside $(t_0 - \ve, t_0 + \ve)$ and $\rhe \to \rh$ uniformly as $\ve \to 0$. 

Let us consider the \Rm\ \mt\ \(g_{\ve}\) on \(\bbr^{n+1}\) defined by \(g_\ve = dr^2 + \rh_\ve(r)^2 \mathrm{g}\), where \(r\) is the Euclidean distance from the origin, and \(\mathrm{g}\) is any metric on the unit sphere \(\mathbb{S}^n = \{r=1\} \subset \bbr^{n+1}\) \w\ $\mathrm{Vol}(\mathbb{S}^n, \mathrm{g}) = 1$. Let us fix $c_0 > n$. We are going to show that if $\ve$ sufficiently small and $c \ge c_0$, then \(W^c(\bbr^{n+1} , g_{\ve}) = 0\). 

For \(R \ge 2\) and \(t \in [0, R]\), let \(\Ps_R(t) = \{t \le r \le R\} \subset \mathbb{R}^{n+1}\), $D_R = \{r < R\}$, $f_{\ve, R}(t) = \mathcal{A}^c_{g_\ve}(\Ps_R(t), D_R)$. Then, for $0 \le t < R$,
\begin{equation}\label{e.f.eps.R}
	\mathcal{A}^c_{g_\ve}(\Ps_R(t), D_R) = f_{\ve, R}(t) = \rhe(t)^n - c \int_{t}^{R} \rhe(s)^n\;ds,
\end{equation}
$f_{\ve, R}(R) = 0$, $\lim\limits_{t \to R^-} f_{\ve, R}(t) = e^{-nR}$. Let
\begin{equation}\label{e.f_R}
	f_{R}(t) = \rh(t)^n - c \int_{t}^{R} \rh(s)^n\;ds,\quad t \in [0, R].
\end{equation}
Then, for $t \neq t_0$, $f'_{R} = \rh^{n-1}(n \rh' + c \rh)$; hence $f_R$ is strictly increasing on $[0, R]$. Further, $f_R(R) = e^{-nR}$. By direct computation, for $1 \le t < R$,
$$ f_{\ve, R}(t) = f_{R}(t) = e^{-nR} - \left( \frac{c}{n} -1 \right)(e^{-nt} - e^{-nR}).$$
So
\begin{equation}\label{e.f_(R)(1)}
	 f_{\ve, R}(1) = f_{R}(1) \le e^{-nR} - \left( \frac{c_0}{n} -1 \right)(e^{-n} - e^{-2n}).
\end{equation}
Since $\rhe \to \rh$ uniformly as $\ve \to 0$ and $\rhe = \rh$ on $[1, \infty)$, using \eqref{e.f.eps.R}, \eqref{e.f_R} and \eqref{e.f_(R)(1)} we conclude that if $\ve$ is sufficiently small,
$$ W^c (D_R, g_\ve) \le \sup_{t \in [0, R]} \mathcal{A}^c_{g_\ve}(\Ps_R(t), D_R) = e^{-nR} \to 0 \text{ as } R \to \infty.$$ 
This implies \(W^c(\bbr^{n+1} , g_\ve) = 0\).
\end{rmk}

\subsection*{Proof of Theorem \ref{main theorem}}

\sps\ \(M\) is a complete \Rm\ \mf\ \w\ sublinear volume growth. It was proved in \cite{cl} that \(M\) satisfies property $(\mathbf{P})$ in Definition \ref{d.prop.P}. Let \(h \in L^1(M) \cap C^\infty(M)\) and \(W^h(M) > 0\). Then, by Lemma \ref{l.existence.of.good.sets}, \tes\ \(\Om \subset M\), which is \(h\)-good. Let \(f: M \to [0, \infty)\) be an unbounded, proper Morse \fn\ and \(\{a_i\}_{i=0}^\infty \subset [0, \infty)\) be a strictly increasing \seq\ \st\ \(\Om \Subset \tilde{\Om}:= f^{-1}([0, a_0 ) ) \), and each \(a_i\) is a regular value of \(f\). 

Using the argument in \cite[Section 12.2]{mon}, we construct a \seq\ of closed \Rm\ \mfs\ \(\{M_i\}_{i=1}^\infty\) as follows. Let \(S_i = f^{-1}(a_i) \). \Tes\ \(\ve_i > 0\) \st\ \(t\) is a regular value of \(f\) \fa\ \(t \in [a_i , a_i + 3 \ve_i] \)  and \tes\ a diffeomorphism \(F_i : f^{-1}([a_i , a_i + 3 \ve_i]) \to S_i \times [0, 3\ve_i]\) which maps \(S_i\) to \(S_i \times \{0\}\). Let \(N_i = f^{-1}([0 , a_i + 3 \ve_i])\). Using a bump \fn\ \(: [a_i , a_i + 3 \ve_i] \to [0,1]\), one can construct a \Rm\ \mt\ \(g_i\) on \(N_i\) \st\ \(g_i = g \) (where \(g\) denotes the \Rm\ \mt\ on \(M\)) on $f^{-1}([0 , a_i +  \ve_i])$ and \(g_i = F_i^* (g|_{S_i} + dt^2)\) (the pullback of the product \mt) on $f^{-1}([a_i +  2\ve_i , a_i +  3\ve_i])$. Finally, we define \(M_i\) to be the double of \(N_i\), i.e. \(M_i\) is obtained by attaching 2 copies of \(N_i\), say \(N_i',\; N_i''\), using the identity map of \(\del N_i\). 

Let \(\Om_i\) denote the isometric copy of \(\Om\) in \(N_i'\). We may also regard \(\Om_i\) as a subset of \(M_i\); in this way we note that
\begin{equation}\label{e.tilde.W.Om_i}
	\widetilde{W}(\Om_i) \leq \sup_{t \in [0, a_0]} \ah(\{ f \le t\}).
\end{equation}
Using a cut-off function, one can construct \(h'_i \in C^{\infty}(N_i)\) \st\
\begin{itemize}
	\item \(h_i' = h\) on $f^{-1}([0 , a_i +  \ve_i])$ and \(h_i = 0\) on $f^{-1}([a_i +  2\ve_i , a_i +  3\ve_i])$;
	\item if \(h_i \in C^{\infty}(M_i)\) is defined by \(h_i = h'_i \) on \(N'_i\) and \(h_i = 0 \) on \(N''_i\), then \(\Om_i \subset M_i\) is \(h_i\)-good and  \(V^{|h_i|}(M_i) \le V^{|h|}(M) + 1\).
\end{itemize}
 \tf\ by Theorems \ref{t.a.m.} and \ref{t.regularity}, \tes\ \((V_i, E_i) \in \mathcal{VC}(M_i)\) \st\ 
 \begin{itemize}
     \item $\ah(V_i, E_i) = \widetilde{W}(\Om_i)$, which implies
\begin{equation}\label{e.H^n.Si_i}
	\|V_i\|(M_i) \le \widetilde{W}(\Om_i) + V^{|h_i|}(M_i) \le \sup_{t \in [0, a_0]} \ah(\{ f \le t\}) + V^{|h|}(M) + 1;
\end{equation}

\item \tes\ \(\Tht \Subset \Om\) and \(\et = \et(\Om, \tilde{\Om})\) \st\ if \(\Tht_i\) denotes the isometric copy of \(\Tht\) inside \(\Om_i\), 
\begin{equation}\label{e.Si_i.Tht_i}
	\|V_i\|(\overline{\Tht}_i) \geq \et;
\end{equation} 

\item $\Si_i = \text{spt} \|V_i\|$ is a closed \(h_i\)-PMC hypersurface \w\ optimal regularity.
 \end{itemize}

Let $\Xi_i$ denote the isometric copy of $\{f < a_i\}$ in $N'_i$, \(V_i' = V_i \mres \Xi_i\), $\Si'_i = \Si_i \cap \Xi_i$. We can view $V_i'$ as an integral $n$-varifold and \(\Si'_i\) as a \(h\)-PMC hypersurface in \(\{f < a_i\} \subset M\). \eqref{e.H^n.Si_i}, \eqref{e.Si_i.Tht_i} imply up to a subsequence, \(V_i'\) converges to an integral $n$-varifold $V_\infty$ in $M$ \st\ 
\begin{equation}
\|V_{\infty}\|(M) < \infty, \quad	\|V_{\infty}\|(\overline{\Tht}) \geq \et.
\end{equation}
Corollary \ref{cor:regularity-stable} and Remark \ref{r.lim.min-max.hyp} imply \(\text{spt} \|V_\infty\|\) is a \(h\)-PMC hypersurface \w\ optimal regularity.

\bibliographystyle{amsalpha}
\bibliography{pmc}
\end{document}